\documentclass[12pt]{amsart}
\usepackage{amsbsy}
\usepackage{graphicx,epsfig,subfigure,psfrag}
\begin{document}

\newtheorem{theorem}{Theorem}
\newtheorem{proposition}{Proposition}
\newtheorem{lemma}{Lemma}
\newtheorem{corollary}{Corollary}
\newtheorem{definition}{Definition}
\newtheorem{remark}{Remark}
\newcommand{\tex}{\textstyle}
\numberwithin{equation}{section} \numberwithin{theorem}{section}
\numberwithin{proposition}{section} \numberwithin{lemma}{section}
\numberwithin{corollary}{section}
\numberwithin{definition}{section} \numberwithin{remark}{section}
\newcommand{\ren}{\mathbb{R}^N}
\newcommand{\re}{\mathbb{R}}
\newcommand{\n}{\nabla}
\newcommand{\p}{\partial}
\newcommand{\iy}{\infty}
\newcommand{\pa}{\partial}
\newcommand{\fp}{\noindent}
\newcommand{\ms}{\medskip\vskip-.1cm}
\newcommand{\mpb}{\medskip}
\newcommand{\AAA}{{\bf A}}
\newcommand{\BB}{{\bf B}}
\newcommand{\CC}{{\bf C}}
\newcommand{\DD}{{\bf D}}
\newcommand{\EE}{{\bf E}}
\newcommand{\FF}{{\bf F}}
\newcommand{\GG}{{\bf G}}
\newcommand{\oo}{{\mathbf \omega}}
\newcommand{\Am}{{\bf A}_{2m}}
\newcommand{\CCC}{{\mathbf  C}}
\newcommand{\II}{{\mathrm{Im}}\,}
\newcommand{\RR}{{\mathrm{Re}}\,}
\newcommand{\eee}{{\mathrm  e}}
\newcommand{\LL}{L^2_\rho(\ren)}
\newcommand{\LLL}{L^2_{\rho^*}(\ren)}
\renewcommand{\a}{\alpha}
\renewcommand{\b}{\beta}
\newcommand{\g}{\gamma}
\newcommand{\G}{\Gamma}
\renewcommand{\d}{\delta}
\newcommand{\D}{\Delta}
\newcommand{\e}{\varepsilon}
\newcommand{\var}{\varphi}
\newcommand{\lll}{\l}
\renewcommand{\l}{\lambda}
\renewcommand{\o}{\omega}
\renewcommand{\O}{\Omega}
\newcommand{\s}{\sigma}
\renewcommand{\t}{\tau}
\renewcommand{\th}{\theta}
\newcommand{\z}{\zeta}
\newcommand{\wx}{\widetilde x}
\newcommand{\wt}{\widetilde t}
\newcommand{\noi}{\noindent}
\newcommand{\uu}{{\bf u}}
\newcommand{\xx}{{\bf x}}
\newcommand{\yy}{{\bf y}}
\newcommand{\zz}{{\bf z}}
\newcommand{\aaa}{{\bf a}}
\newcommand{\cc}{{\bf c}}
\newcommand{\jj}{{\bf j}}
\newcommand{\ggg}{{\bf g}}
\newcommand{\UU}{{\bf U}}
\newcommand{\YY}{{\bf Y}}
\newcommand{\HH}{{\bf H}}
\newcommand{\GGG}{{\bf G}}
\newcommand{\VV}{{\bf V}}
\newcommand{\ww}{{\bf w}}
\newcommand{\vv}{{\bf v}}
\newcommand{\hh}{{\bf h}}
\newcommand{\di}{{\rm div}\,}
\newcommand{\ii}{{\rm i}\,}
\newcommand{\inA}{\quad \mbox{in} \quad \ren \times \re_+}
\newcommand{\inB}{\quad \mbox{in} \quad}
\newcommand{\inC}{\quad \mbox{in} \quad \re \times \re_+}
\newcommand{\inD}{\quad \mbox{in} \quad \re}
\newcommand{\forA}{\quad \mbox{for} \quad}
\newcommand{\whereA}{,\quad \mbox{where} \quad}
\newcommand{\asA}{\quad \mbox{as} \quad}
\newcommand{\andA}{\quad \mbox{and} \quad}
\newcommand{\withA}{,\quad \mbox{with} \quad}
\newcommand{\orA}{,\quad \mbox{or} \quad}
\newcommand{\atA}{\quad \mbox{at} \quad}
\newcommand{\onA}{\quad \mbox{on} \quad}
\newcommand{\ef}{\eqref}
\newcommand{\mc}{\mathcal}
\newcommand{\mf}{\mathfrak}

\newcommand{\ssk}{\smallskip}
\newcommand{\LongA}{\quad \Longrightarrow \quad}
\def\com#1{\fbox{\parbox{6in}{\texttt{#1}}}}
\def\N{{\mathbb N}}
\def\A{{\cal A}}
\newcommand{\de}{\,d}
\newcommand{\eps}{\varepsilon}
\newcommand{\be}{\begin{equation}}
\newcommand{\ee}{\end{equation}}
\newcommand{\spt}{{\mbox spt}}
\newcommand{\ind}{{\mbox ind}}
\newcommand{\supp}{{\mbox supp}}
\newcommand{\dip}{\displaystyle}
\newcommand{\prt}{\partial}
\renewcommand{\theequation}{\thesection.\arabic{equation}}
\renewcommand{\baselinestretch}{1.1}
\newcommand{\Dm}{(-\D)^m}

\title
{\bf Branching analysis of a countable family\\ of global
similarity solutions\\ of a fourth-order thin film equation}

\author{P.~\'Alvarez-Caudevilla and Victor~A.~Galaktionov}

\address{Centro di Ricerca Matematica Ennio De Giorgi,
Scuola Normale Superiore, 56100-Pisa, ITALY}
\email{alvcau.pablo@gmail.com}

\address{Department of Mathematical Sciences, University of Bath,
 Bath BA2 7AY, UK}
\email{vag@maths.bath.ac.uk}


\thanks{The first author is supported by the Ministry of Science and Innovation of
Spain under grant CGL2000-00524/BOS and the Postdoctoral
Fellowship--2008-080.}

\date{\today}





\begin{abstract}

The main goal of the paper is to justify that source-type and
 other global-in-time similarity solutions of the
Cauchy problem for the fourth-order {thin film equation}
  \be
   \label{01}
 u_t=-\n \cdot (|u|^n \n \D u) \inB \ren \times \re_+ \whereA n>0,
 \,\,\, N \ge 1,
  \ee
  can be obtained by a continuous deformation (a homotopy path) as
  $n \to 0^+$ by reducing to similarity solutions (given by eigenfunctions of a rescaled
   linear operator $\BB$) of the classic {\em bi-harmonic
  equation}
   \be
    \label{02}
     \tex{
 u_t = - \D^2 u \inB \ren \times \re_+, \,\,\,\mbox{where}
 \,\,\,\BB=-\D^2 +\frac 14\, y \cdot \n+ \frac N4\, I.
   }
   \ee
  This approach leads to  a countable family of various global similarity patterns of
   \ef{01} and  describes  their oscillatory sign-changing behaviour
   by using the known asymptotic properties
  of the fundamental solution of \eqref{02}.  The branching  from
  $n=0^+$ for \eqref{01} requires Hermitian spectral theory for a
  pair $\{\BB,\BB^*\}$ of non-self adjoint operators and leads to
  a number of difficult mathematical problems. These include, as a key
  part, the problem of multiplicity of solutions, which is under 
 particular scrutiny.

\end{abstract}

\maketitle

\section{Introduction: TFEs, connections with classic PDE theory, layout}
 \label{S1}

\subsection{Main models, their applications, and preliminaries}

\noindent
    We study the global-in-time behaviour of  solutions of
the fourth-order quasilinear evolution equation of parabolic type,
called the {\em thin film equation} (the TFE--4, in short)
\begin{equation}
\label{i1}
    u_{t} = -\nabla \cdot(|u|^{n} \nabla \D u)
 \quad \mbox{in} \quad \ren \times \re_+
    \,,
\end{equation}
where $\n={\rm grad}_x$, $\D=\n \cdot \n$ stands for the Laplace
operator in $\ren$, and $n>0$ is a real parameter. Fourth- and
sixth-order TFEs (the TFE--6) having a similar form,
\begin{equation}
\label{i1166}
    u_{t} = \nabla \cdot (|u|^{n} \nabla \D^2 u)\,,
\end{equation}
 as well as more complicated doubly nonlinear degenerate parabolic models  (see typical examples in \cite{Ki01}),
 have various applications in thin film,
lubrication theory, and in several other hydrodynamic-type
problems. We refer e.g., to \cite{EGK1, EGK3, Gia08, Grun04} for
most recent  surveys and for extended lists of references
concerning physical derivations of various models, key
mathematical results, and further applications. Since the 1980s,
such equations also play a quite special role in nonlinear PDE
theory to be discussed in greater detail below.

 The TFE--4 \eqref{i1} is written for solutions of changing
sign, which can occur in the \emph{Cauchy problem} (the CP) and
also in some \emph{free--boundary problems} (FBPs); see  proper
settings shortly. It is worth mentioning that {\em nonnegative
solutions} with compact support  of various FBPs are mostly
physically relevant, and that the pioneering mathematical
approaches by Bernis and Friedman in 1990 \cite{BF1} were
developed mainly for such solutions.

 However, {\em solutions of
changing sign} have been already under scrutiny for a few years
(see \cite{BW06, EGK2, EGK4}), which in particular can have some
biological motivations \cite{KingPers}, to say nothing of general
PDE theory. It turned out that  these classes of the so-called
``oscillatory solutions of changing sign" of (\ref{i1}) were
rather difficult to tackle rigorously by standard and classic
methods. Moreover, even their self-similar (i.e., ODE)
representatives can lead to several surprises in trying to
describe sign-changing features close to interfaces; see
\cite{EGK2} for a collection of such hard properties. It  turned
out also that, for better understanding of such singular
oscillatory properties of solutions of the CP for (\ref{i1}), it
is  fruitful to consider the (homotopic) limit $n \to 0^+$,
 owing to Hermitian spectral theory developed in
\cite{EGKP} for a pair $\{\BB,\BB^*\}$ of linear rescaled
operators  for $n=0$, i.e., for the {\em bi-harmonic equation}
\begin{equation}
\label{s1}
 \tex{
    u_{t} = -\D^2 u\quad \hbox{in} \quad \ren \times \re_+\,,
    \quad \BB=-\D^2 +\frac 14\, y \cdot \n+ \frac N4\, I, \quad
    \BB^*=-\D^2 - \frac 14 \, y \cdot \n,
    }
\end{equation}
which will always be key for our further analysis.

In the present paper, using this continuity/homotopy deformation
approach ``$n \to 0^+$",
 we shall focus our
analysis to the Cauchy problem for (\ref{i1}) for exponents $n>0$,
which are assumed to be sufficiently small. Some  key and
necessary references will be presented later on.
 We study the large time behaviour of the solutions of \eqref{i1}.
 To this end, we will use some natural connections with a similar
 analysis for more complicated models such as
the {\em limit unstable fourth-order thin film equation} (the
unstable TFE--4):
\begin{equation}
\label{i2}
    u_{t} = -\nabla \cdot(|u|^{n} \nabla \D u)-\D(|u|^{p-1}u)\,,
\end{equation}
with the unstable homogeneous second-order diffusion term, where
$p>1$ is a fixed exponent; see \cite{EGK1} for physical
motivations, references,  and other basics. Here, \eqref{i2}
represents a fourth-order nonlinear parabolic equation with the
backward (unstable) diffusion term in the second-order operator.
Blow-up and global self-similar solutions of \eqref{i2} have been
extensively studied in \cite{EGK1, EGK2} for the unstable TFE--4
(\ref{i2})  and in \cite{EGK3, EGK4} for the {\em unstable
TFE--6},
\begin{equation}
\label{i11}
    u_{t} = \nabla \cdot (|u|^{n} \nabla \D^2 u)-\D(|u|^{p-1}u)\,,
\end{equation}
where further references and other related higher-order TFEs  can
be found.
\par
From the application point of view, it is well known (see
references to surveys above) that \eqref{i1} and \eqref{i2} arises
in numerous areas. In particular, those equations model the
dynamics of a thin film of viscous fluid, as the spreading of a
liquid film along a surface, where $u$ stands the height of the
film (then clearly $u \ge 0$ that naturally leads to a FBP
setting). In particular, when $n=3$ we are dealing with a problem
in the context of lubrication theory for thin viscous films that
are driven by surface tension and when $n=1$ with Hele--Shaw
flows. It is also important to note that, in \ef{i2}, the
fourth-order term reflects surface tension effects and the
second-order term can reflect gravity, van der Waals interactions,
thermocapillary effects, or geometry of the solid substrate.
\par

 Finally, in order to summarize let us mention again that
 higher-order semilinear and quasilinear
parabolic equations occur in applications to thin film theory,
nonlinear diffusion, lubrication theory, flame and wave
propagation (the Kuramoto--Sivashinsky equation and the extended
Fisher--Kolmogorov equation), phase transition at critical
Lifshitz points and bi-stable systems (see Peletier--Troy
\cite{PT} for further details, models, and results). Moreover, in
the special situation when $n=0$ we should notice that \eqref{i2}
is the well known \emph{unstable Cahn--Hilliard equation} (the
CHE)
\begin{equation}
\label{i8}
    u_{t} = - \D^2 u -\D(|u|^{p-1}u) \inB \ren \times \re_+\,;
\end{equation}
 see main references and  full details in \cite{EGW}. The unstable semilinear model
\ef{i8} and similar stable ones are much better known and are
connected with several new applications that have increased the
interest of the study of their quasilinear TFE extensions
\eqref{i2} and consequently of \eqref{i1}.
 Note that, without any doubts, the semilinear CHE \ef{i8} in the
 CP setting,
 admits {\em oscillatory} solutions of changing sign, though with
 no finite (i.e.,
 {\em infinite}) interfaces. As our main goal, we plan to
 extend those properties of the CHE \ef{i8} to the TFEs for small
 $n>0$, where oscillations begin to concentrate at {\em finite} interfaces.


\subsection{A digression to reaction-diffusion theory}

Furthermore,  in the CP setting  for \eqref{i8}, one can write
\eqref{i8} in the form
\begin{equation}
\label{i9}
    \mathcal{A} u_{t} = \D u +u^{p} \whereA  \mathcal{A}:= (-\D)^{-1}
\end{equation}
is a standard  positive operator,
so that \ef{i9} is  a pseudo-parabolic second-order equation. Such
equations had been widely studied since the 1970's, and  nowadays
are well-understood with both the existence and uniqueness results
of local and global classical or blow-up solutions obtained.
first blow-up results for such pseudo-parabolic PDEs were due to
Levine in 1973,  \cite{L}. As was noted in \cite{EGW}, there are
some similarities between \eqref{i9} and the classical semilinear
heat equation from combustion theory (the solid fuel model)
\begin{equation}
\label{i10}
    u_{t} = \D u +u^{p} \quad \hbox{in} \quad \ren \times \re_+,\quad p > 1\,,
\end{equation}
with $N \geq 1$, especially for the blow-up singularity formation
phenomena. Mathematical literature devoted to the study of
\eqref{i10} include a  huge number of papers published since
Fujita's classic papers in the 1960s, and this remarkable history
has been already explained in a dozen of well-known monographs; we
refer to
\cite{BebEb, GalGeom,
AMGV, MitPoh,  Pao, QSupl, SGKM, Sp2}, where further extensions
and references can be found.
In particular, concerning blow-up patterns, a complete description
of all possible types of blow-up have been achieved for some
ranges of the parameters $p$ and $N$, especially in the
subcritical Sobolev range
 $$
 \tex{
 p<p_{\rm Sobolev}= \frac{N+2}{(N-2)_+}.
  }
  $$
 Note that, for $p \ge p_{\rm Sobolev}$, such a
 classification of blow-up patterns is far away from being
 complete, with a number difficult open problems posed.

 Nevertheless, using a standard {\em pseudo-parabolic} form
 \ef{i9} of the C--H flows could rise a hope to apply a huge experience
 achieved earlier for classic reaction-diffusion models \ef{i10}, though
 this is not expected to be that straightforward.


\subsection{A digression to porous medium equation: homotopy to
the heat equation}

Returning to the TFEs, note that the unstable nonlinear operator
in \eqref{i2} gives us the notorious classic \emph{porous medium
equation} (the PME--2, but posed backwards in time),
\begin{equation}
 \label{PME2}
    u_{t} = \D (|u|^{n-1}u)
 \quad(\mbox{for convenience, $p$ is replaced by $n$ as in (\ref{i1})})
    \, ,
\end{equation}
which  derives its name from the role in the description of flows
in porous media.  Parabolic PDE models
  in  filtration theory of liquids and gases in
  porous media were  derived by Leibenzon in the 1920s and 1930s,
 as well as by 
  Richard's  (1931),
   and Muskat (1937).
In fact, modern filtration theory goes back to the beginning of
the twentieth century initiated in the  works by N.Ye.~Zhukovskii,
who is better known for his fundamental research in aerodynamics,
hydrodynamics, and ODE theory (on his pioneering non-oscillation
test in 1892, see \cite[p.~19]{GalGeom}).
   His contribution to ``theory of ground waters"
  is explained in   P.Ya.~Kochina's
 paper \cite{Koch97}.
 For an extended list of references on this subject and more
 filtration history, see
   \cite{GalRDE4n}.

 It is well understood that, for any
$n>1$, \ef{PME2}  has a family of exact self-similar compactly
supported source type solutions (the ZKB ones from 1950s), which
describe the large time behaviour of compactly supported solutions
with conservation of mass in the case on a non-zero mass, i.e.,
\begin{equation*}
 \tex{
    \int\limits u(x,t)\, {\mathrm d}x \neq 0.
    }
\end{equation*}
On the other hand, \ef{PME2} also admits a countable (at least)
family of other similarity solutions; see \cite{GHCo} for key
references and most recent results.

  The PME--2 \ef{PME2}
 can be interpreted as a nonlinear degenerate
version of the classic \emph{heat equation} for $n=0$,
  $$
   u_t=\D u \inA.
   $$
 Note that passing to the limit $n \to 0^+$ in \ef{PME2} for nonnegative solutions
  used to be a difficult
 mathematical problem in the 1970s-80s, which exhibited typical
 (but clearly simpler than in the TFE case) features
 of a ``homotopy" transformation of PDEs.
 This study in 1D was initiated by Kalashnikov in 1978
 \cite{Kal78}. Further detailed results in $\ren$ were obtained in
 \cite{Ben81}; see also \cite{Cock99}. More recent involved
 estimates were obtained in \cite{Pan08, Pan07} for the 1D PME--2 \ef{PME2}
 establishing the rate of convergence of solutions
 as $n \to 0^\pm$, such as $O(n)$ as $n \to 0^-$ (i.e, from $n<0$, the fast diffusion
range, where solutions are smoother)
  in $L^1(\re)$ \cite{Pan08}, and $O(n^2)$ as $n \to 0^+$ in
 $L^2(\re \times (0,T))$ \cite{Pan07}. However, the most of such convergence
 results are obtained for {\em nonnegative} solutions of
 \ef{PME2}.
 For solutions of changing sign, there are some open problems; see \cite{GHCo} for
 references and further details.

Note that, as customary,  any kind of detailed asymptotic analysis
for higher-order equations is much more difficult than for
second-order counterparts in view of the lack of Maximum
Principle, comparison, and order preserving semigroups and
potential properties of the operators involved. Thus, practically
all the existing methods for the
 PME--2 \ef{PME2}
  are not applicable to the TFE--4 \ef{i1}
or \eqref{i2}.
\par

Thus,  in the twenty-first century, higher-order TFEs such as
\ef{i1} and \ef{i1166}, though looking
  like a natural counterpart/extension of the PME--2
 \ef{PME2}, corresponding mathematical
TFE theory is more complicated with several problems,  remaining
open still.

\subsection{Main approaches,  results, and layout}

It is worth mentioning that, unlike the FBPs, studied in hundreds
of papers since 1980s (see \cite{Gia08} and \cite{EGK2} for key
references and  alternative versions of uniqueness approaches),
thin film theory for the Cauchy problem for \eqref{i1} or
\eqref{i2} led recently to a number of difficult open problems and
is not still fully developed; see the above references as a guide
to main difficulties and ideas.
 In fact, the concept of  proper solutions is still rather obscure for the
 Cauchy problem, since any classic or standard notions of
 weak-mild-generalized-... solutions fail in the CP setting.

In this work, we perform a more systematic than before analysis
 of the behaviour of the similarity solutions
through a so-called {\em homotopic approach} (branching from
$n=0$)
via branching theory, using the Lyapunov--Schmidt methods for
obtaining relevant results and properties for the solutions of the
self-similar equation associated with \eqref{i1} and, hence, for
the proper solutions of \eqref{i1}. Overall, loosely speaking,
this approach is characterized as follows: good proper (similarity
or not) solutions of the Cauchy problem for the TFE \eqref{i1} are
those that can be continuously deformed (via a homotopic path) as
$n \to 0^+$ to the corresponding solutions of the {bi-harmonic
equation} \ef{s1},
which will play a crucial role in the subsequent analysis. This
homotopic-like approach
is based upon the spectral properties known for the linear
counterpart \eqref{s1} of the TFE \eqref{i1}. Moreover, owing to
the oscillatory character of the solutions of the bi-harmonic
equation (\ref{s1}) being
a ``limit case" of the TFE \eqref{i1}, close to the interfaces,
this homotopy study exhibits a typical difficulty concerning the
desired structure of the transversal zeros of solutions, at least
for small $n>0$. Proving such a transversality zero property is a
difficult open problem, though qualitatively, this was rather well
understood, \cite{EGK1}.

\par


Indeed, we ascertain through an analytic ``homotopy", which is
understood as just the existence of a continuity as $n \to 0^+$,
over the CP performed in the last section of this paper that the
solutions of \eqref{i1} are homotopic to the solutions of the
bi-harmonic equation \eqref{s1} in a weak sense.
 However, we must admit that this does not  solve the
 problem of uniqueness of solutions of the CP (see details in Section
 \ref{SHom7}), since the final identification of the solutions obtained via analytic
  $\e$-regularization and passing to the limit $\e \to 0^+$ remains not fully understood.
 Overall, it seems that $\e$-regularizations of solutions of the
 TFEs via families of uniformly parabolic analytic flows, which
 was a powerful and successful  tool for second-order degenerate parabolic PDEs
 (such as the PME--2 \ef{PME2}), for TFEs again leads to difficult
 {\em boundary layer}-type problems that remain open in a
 sufficient generality.

\par

 Some parts of the study of the thin film equation \eqref{i2}
 can be performed in similar lines, though a full homotopy
 approach would include the passage to the limit $p \to 1^+$,
 leading to the limit linear equation (not treated here)
 $$
 u_t= -\D^2 u - \D u.
  $$


Thus, the layout of the paper is as follows:


 \noi {\bf (I)} Study of a countable family of global self-similar solutions of
(\ref{i1}) via their branching from eigenspaces at $n=0^+$,
Sections \ref{S2}--\ref{S5}, and


 \noi {\bf (II)} Some general aspects of the CP for (\ref{i1})
by another homotopy approach, Section \ref{SHom7}.

\setcounter{equation}{0}
\section{Problem setting and self-similar solutions}
\label{S2}

\subsection{The FBP and CP}

For  both the FBP and the CP, the solutions are assumed to satisfy
standard free-boundary conditions:
\begin{equation}
\label{i3}
    \left\{\begin{array}{ll} u=0, & \hbox{zero-height,}\\
    \nabla u=0, & \hbox{zero contact angle,}\\
    -{\bf n} \cdot \nabla (|u|^{n}  \D u)=0, &
    \hbox{conservation of mass (zero-flux)}\end{array} \right.
\end{equation}
at the singularity surface (interface) $\Gamma_0[u]$, which is the
lateral boundary of
\begin{equation*}
    \hbox{supp} \;u \subset \ren \times \re_+,\quad N \geq 1\,,
\end{equation*}
where ${\bf n}$ stands for the unit outward normal to
$\Gamma_0[u]$.  Note that, for sufficiently smooth interfaces,
the condition on the flux can be read as
\begin{equation*}
    \lim_{\hbox{dist}(x,\Gamma_0[u])\downarrow 0}
    -{\bf n} \cdot \nabla (|u|^{n}  \D u)=0.
\end{equation*}
It is key that, for the CP, the solutions are assumed to be
``smoother" at the interface than those for the FBP, i.e., \ef{i3}
are not sufficient to define their regularity. These {\em maximal
regularity} issues for the CP, leading to oscillatory solutions,
are under scrutiny in \cite{EGK2}.

Next,  denote by
\begin{equation*}
 \tex{
    M(t):=\int
      u(x,t) \, {\mathrm d}x
    }
\end{equation*}
the mass of the solution, where integration is performed over
smooth support ($\ren$ is allowed for the CP only). Then,
differentiating $M(t)$ with respect to $t$ and applying the
Divergence Theorem (under natural regularity assumptions on
solutions and free boundary), we have that
\begin{equation*}
 \tex{
  J(t):=  \frac{{\mathrm d}M}{{\mathrm d}t}= -
  \int\limits_{\Gamma_0\cap\{t\}}{\bf n} \cdot \nabla
     (|u|^{n}  \D u )\, .
     }
\end{equation*}
The mass is conserved if
   $ J(t) \equiv 0$, which is assured by the flux condition in
   \eqref{i3}.
\par
The problem is completed with bounded, smooth, integrable,
compactly supported initial data
\begin{equation}
\label{i4}
    u(x,0)=u_0(x) \quad \hbox{in} \quad \Gamma_0[u] \cap \{t=0\}.
\end{equation}

In the CP for \eqref{i1} in $\ren \times \re_+$, one needs to pose
bounded compactly supported initial data \eqref{i4} prescribed in
$\ren$. Then,  under the same zero flux condition at finite
interfaces (to be established separately), the mass is preserved.
\subsection{Global similarity solutions: towards a nonlinear eigenvalue problem}

 We now begin to specify the self-similar
solutions
of the equation
\eqref{i1},
which
are admitted due to its natural scaling-invariant  nature.
In the case of the mass being conserved, we have global in time
source-type solutions.
\par
Using the following scaling in \eqref{i1}
\begin{equation*}
    x:= \mu \bar x,\quad t:= \l \bar t,\quad u:= \nu \bar u, \quad
    \mbox{with}
\end{equation*}
\begin{equation}
\label{sf1}
 \tex{
    \frac{\p u}{\p t}= \frac{\nu}{\l} \frac{\p \bar u}{\p t},\quad
    \frac{\p u}{\p x_{i}}= \frac{\nu}{\mu} \frac{\p \bar u}{\p x_{i}},\quad
    \frac{\p^2 u}{\p x_{i}^2}= \frac{\nu}{\mu^2} \frac{\p^2 \bar u}{\p
    x_{i}^2},
     }
\end{equation}
 and substituting those expressions in \eqref{i1} yields
\begin{equation*}
 \tex{
    \frac{\nu}{\l} \frac{\p \bar u}{\p t}=-
    \frac{\nu^{n+1}}{\mu^4} \nabla \cdot (|\bar u|^{n} \nabla \D \bar u)\,.
    }
\end{equation*}
To keep this equation invariant, the following  must be fulfilled:
\begin{equation}
\label{sf2}
 \tex{
    \frac{\nu}{\l}=\frac{\nu^{n+1}}{\mu^4}, \quad \mbox{so that}
    }
\end{equation}
\begin{equation*}
    \mu := \l^\b \Longrightarrow \nu := \l^{\frac{4\b-1}{n}} \quad
    \mbox{and}\quad
 \tex{
    u(x,t):= \l^{\frac{4\b-1}{n}} \bar u(\bar x,\bar t) = \l^{\frac{4\b-1}{n}}
    \bar u(\frac{x}{\mu},\frac{t}{\l}).
    }
\end{equation*}
Consequently,
\begin{equation*}
 \tex{
    u(x,t):= t^{\frac{4\b-1}{n}} f(\frac{x}{t^\b}),
     }
\end{equation*}
where $t=\l$ and $f(\frac{x}{t^\b})=\bar u(\frac{x}{t^\b},1)$.
Owing to \eqref{sf2}, we obtain
\begin{equation*}
   n\a +4\b=1,
\end{equation*}
which links the parameters $\a$ and $\b$. Hence, substituting
\begin{equation}
\label{sf3}
 \tex{
    u(x,t):= t^{-\a} f(y), \quad \hbox{with}\quad
    y:=\frac{x}{t^\b},\quad \b= \frac{1-n \a}4,
    }
\end{equation}
into \eqref{i1} and rearranging terms, we find that the function
$f$ solves a quasilinear elliptic equation of the form
\begin{equation}
\label{sf4}
    \nabla \cdot(|f|^{n} \nabla \D f)=\a f+\b y \nabla \cdot f\,.
\end{equation}

 Finally,
 thanks to the above relation between $\a$ and $\b$, we find a {\em
nonlinear eigenvalue problem} of the form
\begin{equation}
\label{sf5}
 \fbox{$
 \tex{
    -\nabla \cdot(|f|^{n} \nabla \D f) +\frac{1-\a n}{4}\, y \nabla \cdot f +\a
    f=0,
    \quad f \in C_0(\ren)\, ,
    }
    $}
\end{equation}
 where we add to the equation \ef{sf4} a natural assumption that
 $f$ must be compactly supported (and, of course, sufficiently
 smooth at the interface, which is an accompanying question to be
 discussed as well).

 Thus, for such degenerate elliptic equations,
  the functional setting in \ef{sf5} assumes that we are
 looking for  (weak) {\em compactly supported} solutions $f(y)$ as
 certain ``nonlinear eigenfunctions" that hopefully occur for special values of nonlinear eigenvalues
  $\{\a_k\}_{k \ge 0}$. Our goal is to justify
  that, labelling  the eigenfunctions via a multiindex $\s$,
 \be
 \label{sf51}
 \fbox{$
 \mbox{
 (\ref{sf5}) possesses  a countable set of
 eigenfunction/value pairs $\{f_\s,\, \a_k\}_{|\s|=k \ge 0}$.
  }
  $}
  \ee

 Concerning  the well-known properties of finite propagation for TFEs, we refer to papers
 \cite{EGK1}--\cite{EGK4}, where a large amount of earlier
 references is available; see also \cite{GMPSobI} for more recent
 results and references in this elliptic area. However,  there are still
a little of entirely rigorous results, especially those that are
attributed to the Cauchy problem for TFEs. In the linear case
$n=0$,
 the condition $f \in C_0(\ren)$, is naturally replaced by the requirement that the
 eigenfunctions $\psi_\b(y)$ exhibit typical exponential decay at
 infinity, a property that is reinforced by introducing  appropriate weighted $L^2$-spaces.
Actually,
 using the homotopy limit $n \to 0^+$, we will be obliged for
 small $n>0$,
 instead of $C_0$-setting in
\eqref{sf5}, we will use the following weighted $L^2$-space:
  \begin{equation}
   \label{WW11}
  f \in L^2_\rho(\ren), \quad \mbox{where} \quad \rho(y)={\mathrm e}^{a |y|^{4/3}},
  \quad a>0 \,\,\,\mbox{small}.
  \end{equation}

Note that, in the case of the Cauchy problem with conservation of
mass making use of the self-similar solutions \eqref{sf3}, we have
that
\begin{equation*}
 \tex{
    M(t):=\int\limits_{\ren} u(x,t) \, {\mathrm d}x=t^{-\a} \int\limits_{\ren}
     f\big(\frac{x}{t^\b}\big)
    \, {\mathrm d}x = t^{-\a+\b N} \int\limits_{\ren} f(y)
    \, {\mathrm d}y,
    }
\end{equation*}
 where the actual integration is performed over the support
  ${\rm supp}\, f$
of the nonlinear eigenfunction.
 Then, as is well known,
if $\int f \not = 0$,
  the exponents are calculated giving the first explicit nonlinear eigenvalue:
\begin{equation}
 \label{alb1}
 \tex{
 -\a + \b N=0 \LongA
    \a_0(n)=\frac{N}{4+Nn} \quad \mbox{and}  \quad \b_0(n)=\frac{1}{4+Nn}.
    }
\end{equation}

\setcounter{equation}{0}
\section{Hermitian spectral theory of the linear rescaled operators}
\label{S3}

 In this section, we establish the spectrum $\s(\bf{B})$
of the linear operator $\bf{B}$ obtained from the rescaling of the
linear counterpart of \eqref{i1}, the bi-harmonic equation
(\ref{s1}), which will be  essentially used in what follows.

\subsection{How the operator $\BB$ appears: a linear eigenvalue problem}

Let $u(x,t)$ be the unique solution of the CP for the linear
parabolic bi-harmonic equation \eqref{s1} with the initial data
(the space as in \eqref{WW11} to be more properly introduced
shortly)
\begin{equation*}
    u_0 \in L_{\rho}^2(\ren),
\end{equation*}
given by the convolution Poisson-type integral
\begin{equation}
\label{s2}
 \tex{
    u(x,t)=b(t)\, * \, u_0 \equiv t^{-\frac N4} \int\limits_{\ren} F((x-z)t^{-\frac 14})
     u_0(z)\, {\mathrm d}z.
    }
\end{equation}
Here, by scaling invariance of the problem, in a similar way as
was done in the previous section for \eqref{i1}, the  unique the
fundamental solution of the operator $\frac{\p}{\p t} + \D^2$ has
the self-similar structure
\begin{equation}
\label{s3}
 \tex{
    b(x,t)=t^{-\frac N 4} F(y), \quad y:=\frac{x}{t^{1/4}} \quad  (x\in \ren).
    }
\end{equation}
Substituting $b(x,t)$ into \eqref{s1} yields that the rescaled
fundamental kernel $F$ in \ef{s3} solves the linear elliptic
problem
\begin{equation}
\label{s4}
 \tex{
    {\bf B}F \equiv -\D_y^2 F + \frac{1}{4}\, y \cdot \nabla_y F  +\frac{N}{4}\, F=0
    \quad \hbox{in} \quad \ren,\quad \int\limits_{\ren} F(y) \, {\mathrm
    d}y=1.
    }
\end{equation}

${\bf B}$ is a non-symmetric linear operator, which is bounded
from $H_{\rho}^4(\ren)$ to $L_{\rho}^2(\ren)$ with the exponential
weight as in \ef{WW11}.
Here,  $a\in (0,2d)$ is any positive constant, depending on the
parameter $d>0$, which characterises  the exponential decay of the
kernel $F(y)$:
 \begin{equation}
  \label{F11}
 |F(y)| \le D \, {\mathrm e}^{-d |y|^{4/3}} \quad \mbox{in} \quad
 \ren \quad \big(D>0, \,\,\, d=3 \cdot 2^{-\frac{11}3}\big).
   \end{equation}
 Later on, by  $F$ we denote the oscillatory rescaled kernel as
the only solution of \eqref{s4}, which has exponential decay,
oscillates as $|y|\rightarrow \infty$, and satisfies the  standard
pointwise estimate \ef{F11}.

\ssk

Thus, we need to solve the corresponding {\em linear eigenvalue
problem}:
 \be
 \label{LP1}
  \fbox{$
  \BB \psi = \l \psi \inB \ren, \quad \psi \in L^2_\rho(\ren).
   $}
   \ee
   One can see that the nonlinear
   one (\ref{sf5}) formally reduces to (\ref{LP1}) at $n=0$ with
   the following shifting of the corresponding eigenvalues:
    $$
    \tex{
    \l=-\a + \frac N4.
    }
    $$
 In fact, this is the main reason to calling (\ref{sf5}) a
 nonlinear eigenvalue problem, and, crucially, the discreteness of
 the real spectrum of the linear one (\ref{LP1}) will be shown to
 be inherited by the nonlinear problem, but  a long way is needed
 to justify such an issue.

\subsection{Functional setting and semigroup expansion}

    Thus,
 we solve (\ref{LP1}) and calculate the spectrum of $\s({\bf B})$ in the
weighted space $L_{\rho}^2(\ren)$. We then need the following
Hilbert space:
\begin{equation*}
    H_{\rho}^4(\ren) \subset L_{\rho}^2(\ren) \subset L^2(\ren).
\end{equation*}
The Hilbert space $H_{\rho}^4(\ren)$ has the following inner product:
\begin{equation*}
 \tex{
    \big\langle v,w \big\rangle_{\rho} := \int\limits_{\ren} \rho(y) \sum\limits_{k=0}^{4}
     D^{k} v(y) \overline{{D^{k} w(y)}}
     \, {\mathrm d}y,
      }
\end{equation*}
where $D^{k} v$ stands for the vector
    $\{D^{\b} v\,,\,|\b|=k\}$,
and the norm
\begin{equation*}
 \tex{
    \| v \|_{\rho}^2 := \int\limits_{\ren} \rho(y) \sum\limits_{k=0}^{4}
     |D^{k} v(y)|^2 \, {\mathrm d}y.
    }
\end{equation*}

\par
 Next, introducing the rescaled variables
\begin{equation}
\label{s6}
 \tex{
    u(x,t)=t^{-\frac N4} w(y,\tau), \quad y:=\frac{x}{t^{1/4}}, \quad \tau= \ln t \,:\,
    \re_+ \to \re,
  }
\end{equation}
we find that the rescaled solution $w(y,\t)$ satisfies the
evolution equation
\begin{equation}
\label{s7}
    w_{\tau} = {\bf B}w\,,
\end{equation}
since, substituting the representation of $u(x,t)$ \eqref{s6} into
\eqref{s1} yields
\begin{equation*}
 \tex{
    -\D_y^2 w + \frac{1}{4} \, y \cdot \nabla_y w  +\frac{N}{4} \, w= t\, \frac{\p w}{\p t} \frac{\p \tau}{\p t}.
  }
\end{equation*}
Thus, to keep this invariant, the following should be satisfied:
\begin{equation*}
 \tex{
    t \, \frac{\p \tau}{\p t}=1 \,\, \Longrightarrow \,\, \tau = \ln t, \quad \mbox{i.e.,
    as defined in (\ref{s6})}.
    }
\end{equation*}
Hence, $w(y,\tau)$ is the solution of the Cauchy problem for the
equation \eqref{s7} and with the following initial condition at
$\tau=0$, i.e., at $t=1$:
\begin{equation}
\label{s8}
    w_0(y) = u(y,1)\equiv b(1)\, * \, u_0 = F\, * \, u_0\, .
\end{equation}
Thus, the linear operator $\frac{\p}{\p \tau} - {\bf B}$ is also a
rescaled version of the standard parabolic one $\frac{\p}{\p t} +
\D^2$. Therefore, the corresponding semigroup ${\mathrm e}^{{\bf
B} \tau}$ admits an explicit integral representation. This helps
to establish some properties of the operator ${\bf B}$ and
describes other evolution features of the linear flow. From
\eqref{s2} we find the following explicit representation of the
semigroup:
\begin{equation}
\label{s9}
 \tex{
    w(y,\tau)=\int\limits_{\ren} F \big(y-z{\mathrm e}^{-\frac{\tau}{4}}\big)\, u_0(z) \,
    {\mathrm d}z \equiv {\mathrm e}^{{\bf B} \tau} w_0, \quad \mbox{where}
     \quad
    x=t^{\frac{1}{4}}y,  \quad \tau=\ln t.
    }
\end{equation}
Subsequently, consider Taylor's power series of the analytic
kernel\footnote{We hope that returning here to the multiindex $\b$
instead of $\s$ in \ef{sf51} will not lead to a confusion with the
exponent $\b$ in self-similar scaling \ef{sf3}.}
\begin{equation}
\label{s10}
 \tex{
    F\big(y-z {\mathrm e}^{-\frac{\tau}{4}}\big)=\sum\limits_{(\b)} {\mathrm e}^{
    -\frac{|\b|\tau}{4}} \frac{(-1)^{|\b|}}{\b!}    D^\b F(y) z^\b
    \equiv \sum\limits_{(\b)} {\mathrm e}^{-\frac{|\b|\tau}{4}} \frac{1}{\sqrt{\b!}} \psi_\b(y) z^\b,
    }
\end{equation}
for any $y\in \ren$, where
\begin{equation*}
    z^{\b}:=z_1^{\b_1}\cdots z_{N}^{\b_{N}}
\end{equation*}
and $\psi_{\b}$ are the normalized eigenfunctions for the operator
$\bf{B}$.
The series in \ef{s10} converges uniformly on compact subsets in
$z \in \ren$. Indeed, denoting $|\b|=l$ and estimating the
coefficients
\begin{equation*}
 \tex{
    \big|\sum\limits_{|\b|=l} \frac{(-1)^{l}}{\b!}  D^\b F(y) z_1^{\b_1}\cdots z_{N}^{\b_{N}}\big| \leq b_l
    |z|^l,
    }
\end{equation*}
by Stirling's formula we have that, for $l\gg 1$,
\begin{equation}
\label{s11}
 \tex{
    b_l = \frac{N^l}{l!} \sup_{y\in \ren,
    |\b|=l} |D^\b F(y)| \approx \frac{N^l}{l!} l^{-l/4}
    {\mathrm e}^{l/4} \approx l^{-3l/4}c^l ={\mathrm e}^{-l \ln 3l/4 +l \ln c}.
    }
\end{equation}
Note that, the series
\begin{equation*}
 \tex{
    \sum b_l |z|^l
    }
\end{equation*}
has the radius of convergence  $R=\infty$.

 Thus, we obtain the
following representation of the solution:
\begin{equation}
\label{s12}
 \tex{
    w(y,\tau)= \sum\limits_{(\b)}  {\mathrm e}^{-\frac{|\b|}{4}\, \t} M_\b(u_0) \psi_\b(y),
    \quad \mbox{where} \quad
    \l_\b =: -\frac{|\b|}{4}
    }
\end{equation}
 and $\{\psi_\b\}$ are the eigenvalues and
eigenfunctions of the operator ${\bf B}$, respectively, and
\begin{equation*}
 \tex{
    M_\b(u_0):= \frac{1}{\sqrt{\b!}} \int\limits_{\ren} z_1^{\b_1}
    \cdots z_{N}^{\b_{N}} u_0(z) \, {\mathrm d}z
 }
\end{equation*}
are the corresponding momenta of the initial datum $w_0$ defined
by \eqref{s8}.

\subsection{Main spectral properties of the pair $\{\BB,\,
\BB^*\}$}

 Thus,  the
following holds \cite{EGKP}:

\begin{theorem}
\label{Th s1} {\rm (i)} The spectrum of ${\bf B}$ comprises real
eigenvalues only with the form
\begin{equation}
\label{s13}
 \tex{
    \s({\bf B}):=\big\{\l_\b = -\frac{|\b|}{4}\,,\,|\b|=0,1,2,...\big\}.
    }
\end{equation}
Eigenvalues $\l_\b$ have finite multiplicity with eigenfunctions,
\begin{equation}
\label{s14}
 \tex{
    \psi_\b(y):= \frac{(-1)^{|\b|}}{\sqrt{\b!}} D^\b F(y) \equiv \frac{(-1)^{|\b|}}{\sqrt{\b!}}
    \big(\frac{\p}{\p y_1}\big)^{\b_1}\cdots \big(\frac{\p}{\p y_N}\big)^{\b_N} F(y).
    }
\end{equation}

\noi{\rm (ii)} The subset of eigenfunctions
  $  \Phi=\{\psi_\b\}$
is complete in $L^2(\ren)$ and in $L_{\rho}^2(\ren)$.

\noi{\rm (iii)} For any $\l \notin \s({\bf B})$, the resolvent
  $  ({\bf B}-\l I)^{-1}$
is a compact operator in $L_{\rho}^2(\ren)$.
\end{theorem}

Subsequently, it was also shown in \cite{EGKP} that the adjoint
(in the dual metric of $L^2(\ren)$) operator of ${\bf B}$ given by
\begin{equation}
 \label{B*}
 \tex{
    {\bf B}^*:=-\D^2 - \frac{1}{4}\,\, y \cdot \nabla,
    }
\end{equation}
 in the weighted space
$L_{{\rho}^*}^2(\ren)$, with the exponentially decaying weight
function
\begin{equation}
 \label{rho*}
 \tex{
    {\rho}^*(y) \equiv \frac{1}{\rho(y)} = {\mathrm e}^{-a|y|^\a} >0,
    }
\end{equation}
is a bounded linear operator,
\begin{equation*}
    {\bf B}^*:H_{{\rho}^*}^4(\ren) \to L_{{\rho}^*}^2(\ren),
    \,\, \mbox{so} \,\,
    \big\langle {\bf B} v, w\big\rangle = \big\langle v, {\bf B}^* w\big\rangle,
    \,\,
    v \in H_{\rho}^4(\ren), \,\,
    w \in H_{{\rho}^*}^4(\ren).
\end{equation*}
Moreover, the following theorem establishes the spectral properties of the
adjoint operator which will be very similar to those shown in Theorem\,\ref{Th s1}
for the operator $\bf{B}$.

\begin{theorem}
\label{Th s2}
{\rm (i)} The spectrum of ${\bf B}^*$ consists of eigenvalues of
finite multiplicity,
\begin{equation}
\label{s15}
 \tex{
    \s({\bf B}^*)=\s({\bf B}):=\big\{\l_\b = -\frac{|\b|}{4}\,,\,|\b|=0,1,2,...\big\},
    }
\end{equation}
and the eigenfunctions
   $\psi_\b^*(y)$
are polynomials of order $|\b|$.

\noi{\rm (ii)}  The subset of eigenfunctions
$    \Phi^*=\{\psi_\b^*\}$
is complete
 in $L_{{\rho}^*}^2(\ren)$.

\noi{\rm (iii)} For any $\l \notin \s({\bf B}^*)$, the resolvent
  $  ({\bf B}^*-\l I)^{-1}$
is a compact operator in $L_{{\rho}^*}^2(\ren)$.
\end{theorem}

It should be pointed out that, since $\psi_0=F$ and $\psi_0^*
\equiv 1$, we have
\begin{equation*}
 \tex{
     \langle \psi_0, \psi_0^* \rangle= \int\limits_{\ren} \psi_0\, {\mathrm
     d}y
      =\int\limits_{\ren} F (y) \, {\mathrm
     d}y=1.
     }
\end{equation*}
However, thanks to \eqref{s14}, we have that
\begin{equation*}
 \tex{
     \int\limits_{\ren} \psi_\b \equiv \langle \psi_\b, \psi_0^* \rangle =0 \quad \hbox{for any} \quad
     |\b|\neq 0.
     }
\end{equation*}
This expresses the orthogonality property to the adjoint
eigenfunctions in terms of the dual inner product.

Note that \cite{EGKP}, for the eigenfunctions $\{\psi_\b\}$ of
$\bf{B}$ denoted by \eqref{s14}, the corresponding adjoint
eigenfunctions are {\em generalized Hermite  polynomials} given by
\begin{equation}
\label{s16}
 \tex{
    \psi_\b^*(y):=\frac{1}{\sqrt{\b!}}\Big[y^\b + \sum\limits_{j=1}^{[|\b|/4]}
    \frac{1}{j!}\, \D^{2j} y^\b\Big].
    }
\end{equation}
Hence, the orthonormality condition holds
\begin{equation}
\label{s17}
    \big\langle \psi_\b,\psi_\g \big\rangle=\d_{\b\g}
    \quad \hbox{for any} \quad \b,\;\g,
\end{equation}
where $\big\langle \cdot,\cdot \big\rangle$ is the duality product
in $L^2(\ren)$ and $\d_{\b\g}$ is  Kronecker's delta. Also,
operators $\bf{B}$ and $\bf{B}^*$ have zero Morse index (no
eigenvalues with positive real parts are available).

 Key spectral results  can be extended \cite{EGKP} to
$2m$th-order linear poly-harmonic flows
 \begin{equation}
  \label{PP1}
  \tex{
 u_t= - (-\D)^m u \quad \mbox{in} \quad \ren \times \re_+,
  }
  \end{equation}
  where
 the elliptic equation for the rescaled kernel $F(y)$ takes the
 form
\begin{equation}
\label{s5}
  \tex{
    {\bf B} F \equiv -(-\D_y)^m F + \frac{1}{2m}\,y \cdot \nabla_y F  +\frac{N}{2m} \,F=0
    \quad \hbox{in} \quad \ren,\quad \int\limits_{\ren} F(y) \, {\mathrm
    d}y=1.
     }
\end{equation}
In particular, for $m=1$, we find the \emph{Hermite operator} and
the {\em Gaussian kernel} (see \cite{BS} for further information)
\begin{equation*}
 \tex{
    {\bf B} F \equiv \D F + \frac{1}{2}\, y\cdot \n F  +\frac{N}{2} \,F=0
     \LongA F(y)= \frac 1{(4 \pi)^{N/2}} \, {\mathrm
     e}^{-\frac{|y|^2}4},
    }
\end{equation*}
whose name is connected with fundamental works of Charles Hermite
on orthogonal polynomials $\{H_\b\}$ about 1870.
 These classic Hermite polynomials are obtained by differentiating
 the Gaussian: up to normalization constants,
 \be
 \label{Ga1}
  \tex{
  D^\b{\mathrm
     e}^{-\frac{|y|^2}4}= H_\b(y) \, {\mathrm
     e}^{-\frac{|y|^2}4} \quad \mbox{for any} \,\,\, \b.
     }
     \ee
 Note that, for
$N=1$, such operators and polynomial eigenfunctions in 1D were
studied earlier by Jacques C.F.~Sturm in 1836; on this history and
Sturm's main original calculations, see \cite[Ch.~1]{GalGeom}.

The generating formula \ef{Ga1} for (generalized) Hermite
polynomials is not available if $m \ge 2$, so that \ef{s16} are
obtained via a different procedure, \cite{EGKP}.

\setcounter{equation}{0}
\section{Similarity profiles for the Cauchy problem via  $n$-branching}
\label{S5}

 \subsection{Derivation of the branching equation}

 In general,   construction of oscillatory similarity solutions of the
Cauchy problem for the TFE--4 \ef{i1} is a  difficult nonlinear
problem, which is harder than for the corresponding  FBP one. On
the other hand, for $n=0$, such similarity profiles  exist and are
given by eigenfunctions $\{\psi_\b\}$. In particular, the first
mass-preserving profile is just the rescaled kernel $F(y)$, so it
is unique, as was shown in Section \ref{S3}. Hence, somehow, there
can be expected a possibility to visualize such an oscillatory
first ``nonlinear eigenfunction" $f(y)$ of changing sign, which
satisfies the {\em nonlinear eigenvalue problem} \eqref{sf5},
  at least, for
sufficiently small $n > 0$. This assumes using the $n$-branching
approach that ``connects" $f$ with the rescaled fundamental
profile $F$ satisfying the corresponding linear equation \ef{s4},
with all the necessary properties of $F$  presented in Section
\ref{S3}.

\par
 Thus, we plan to describe the behaviour of the similarity
 profiles $\{f_\b\}$, as nonlinear eigenfunctions of \ef{sf5}
for the TFE performing a ``homotopic" approach when $n\downarrow
0$.
 Homotopic approaches are well-known in the theory of
vector fields, degree, and nonlinear operator theory (see
\cite{Deim, KZ} for details). However, we shall be less precise in
order to apply that approach, and here a ``homotopic path" just
declares existence of a continuous connection (a curve) of
solutions $f \in C_0$ that ends up at $n=0^+$ at the linear
eigenfunction $\psi_0(y)=F(y)$ or further eigenfunctions
$\psi_\b(y) \sim D^\b F(y)$, as \ef{s14} claims.
\par

Using classical branching theory in the case of finite regularity
of nonlinear operators involved, we formally show that the
necessary orthogonality condition holds deriving the corresponding
{\em Lyapunov--Schmidt branching equation}. We will try to be as
much rigorous as possible in supporting of delivering  the
nonlinear eigenvalues $\{\a_k\}$. Further extensions of solutions
for non-small $n>0$ require a novel essentially non-local
technique of such nonlinear analysis, which remains an open
problem.
\par

Those critical eigenvalues $\{\a_k\}$ are obtained according to
non-self-adjoint spectral theory from Section \ref{S3}. We then
use  the explicit expressions for the eigenvalues and
eigenfunctions of the linear eigenvalue problem \ef{LP1} given in
Theorem \ref{Th s1}, where we also need the main conclusions of
the ``adjoint" Theorem \ref{Th s2}.
  Then, taking the corresponding
linear equation from \eqref{sf5} with $n=0$, we find, at least,
formally, that
\begin{equation*}
 \tex{
  n=0: \quad  \mathcal{L}(\a)f:=-\D^2 f +\frac{1}{4}\, y \cdot\nabla f  +\a f=0.
     }
\end{equation*}

From that equation combined with the eigenvalues expressions
obtained in the previous section, we ascertain the following
critical values for the parameter $\a_k=\a_k(n)$,
\begin{equation}
\label{bf4}
 \tex{
  n=0: \quad   \a_k(0) := -\l_k + \frac{N}{4} \equiv \frac{k+N}4 \quad \hbox{for any} \quad k=0,1,2,\ldots,
     }
\end{equation}
where $\l_k$ are the eigenvalues  defined in Theorem\,\ref{Th s1},
so that
\begin{equation*}
 \tex{
    \a_0(0)= \frac{N}{4}, \; \a_1(0)= \frac{N+1}{4},\; \a_2(0)= \frac{N+2}{4},\ldots ,
    \a_k(0)= \frac{k+N}{4}\ldots \,.
     }
\end{equation*}
In particular, when $k=0$, we have that $\a_0(0)= \frac{N}{4}$ and
the eigenfunction satisfies
\begin{equation*}
    {\bf B} F=0, \quad \mbox{so that} \quad
    \ker \mathcal{L}(\a_0) = \mathrm{span\,}\{\psi_0\} \quad(\psi_0=F),
\end{equation*}
and, hence, since $\l_0=0$ is a simple eigenvalue for the operator
$\mathcal{L}(\a_0)= {\bf B}$, its algebraic multiplicity is 1. In
general, we find that
\begin{equation}
\label{bf5}
 \tex{
    \ker\big({\bf B} + \frac{k}{4}\, I \big)= \mathrm{span\,}\{\psi_\b, \, |\b|=k
    \}, \quad \hbox{for any}
    \quad k=0,1,2,3,\cdots\,,
    }
\end{equation}
where the operator ${\bf B} + \frac{k}{4}\, I$ is Fredholm of
index zero since it is a compact perturbation of the identity of
linear type with respect to $k$. In other words,
$R[\mathcal{L}(\a_k)]$ is a closed subspace of $L_{\rho}^2(\ren)$
and, for each $\a_k$,
\begin{equation*}
    \hbox{dim} \ker(\mathcal{L}(\a_k))< \infty \andA \hbox{codim}R[\mathcal{L}(\a_k)]< \infty.
\end{equation*}

Then, for small $n>0$ in \eqref{sf5}, we can use the asymptotic
expansions
\begin{equation}
\label{br3}
    \a_k(n):= \a_k+ \mu_{1,k} n+ o(n),\quad \mbox{and}
 \ee
 \be
  \label{br3N}
     |f|^n \equiv  {\mathrm e}^{n\ln |f|}:=
     1 +n \ln |f|+o(n).
\end{equation}
 As customary in bifurcation-branching theory \cite{KZ, VainbergTr},
existence of an expansion such as \ef{br3} will allow one to get
further expansion coefficients in
\begin{equation}
\label{br31}
    \a_k(n):= \a_k+ \mu_{1,k} n + \mu_{2,k} n^2+ \mu_{3,k} n^3 + ...\,,
 \ee
 as the regularity of nonlinearities allows and suggests, though the convergence
 of such an analytic series can be questionable and is not under
 scrutiny here.

Another principle question is that, for oscillatory sign changing
profiles $f(y)$, the last expansion \ef{br3N} cannot be understood
in the pointwise sense, but can be naturally expected to be valid
in other metrics such as weighted $L^2$ or Sobolev spaces as in
Section \ref{S3} that used to be appropriate for the functional
setting of the equivalent integral equation and for that with
$n=0$. Since \ef{br3N} is obviously pointwise violated at the
nodal set $\{f=0\}$ of $f(y)$, this imposes some restrictions on
the behaviour   of corresponding eigenfunctions $\psi_\b(y)$
($n=0$) close to their zero sets.
 Using well-known asymptotic and other related properties of the
{\em radial} analytic rescaled kernel $F(y)$ of the fundamental
solutions \ef{s3}, the generating formula of eigenfunctions
\ef{s14} confirms that the nodal set of analytic eigenfunctions
$\{\psi_\b=0\}$
 consists of isolated zero surfaces, which are ``transversal", at least in the a.e. sense,
with the only accumulation point at $y= \infty$.
 Overall, under such conditions, this indicates that
  \be
  \label{log1}
 \mbox{expansion (\ref{br3N}) contains not more than
 ``logarithmic" singularities a.e.},
  \ee
   which well suited the integral compact operators involved into
   a branching analysis,
  though we are far away to claim this as any rigorous issue.
 Moreover, when $n>0$ is not small enough, such an analogy and statements like
  \ef{log1} become not that
clear, and global extensions of continuous $n$-branches induced by
some compact integral operators, i.e., nonexistence of turning
(saddle-node) points in $n$, require, as usual, some unknown
monotonicity-like results.

\par

Now, we assume the expansion \ef{br3N} away from possible zero
surfaces of $f(y)$, which, by transversality, can be localized in
arbitrarily small neighbourhoods. Indeed, it is clear that when
$|f| > \d >0$, for any $\d>0$, there is no problem in
approximating $|f|^n$ by \eqref{br3}, i.e., $|f|^n =1+ O(n)$ as $n
\rightarrow 0^+$. However, when $|f| \leq \d $, for any $\d \geq
0$ sufficiently small, the proof of such an approximation in weak
topology (as suffices for dealing with equivalent integral
equations) is far from clear unless the zeros of the $f$'s are
also transversal a.e. with a standard accumulating property at the
only interface zero surface. The latter issues have been studied
and described in \cite{EGK2} in the radial setting.
 Hence, we can suppose that such nonlinear eigenfunctions $f(y)$ are oscillatory
and infinitely sign changing close to the interface surface.
Therefore, if we assume that their zero surface are transversal
a.e. with a known geometric-like accumulation at the interface, we
find that, for any $n$ close to zero and any $\d= \d(n)
>0$ sufficiently small,
\begin{equation*}
    n| \ln |f| | \gg 1, \quad \hbox{if} \quad |f| \leq \d(n),
\end{equation*}
and, hence, on such subsets, $f(y)$ must be exponentially small:
\begin{equation*}
 \tex{
    | \ln |f| | \gg \frac{1}{n}\; \Longrightarrow \;\ln |f| \ll -\frac{1}{n}\;
    \Longrightarrow \; |f|  \ll {\mathrm e}^{-\frac{1}{n}}.
    }
\end{equation*}
Recall that this happens in  also exponentially small
neighbourhoods of the transversal zero surfaces.

Overall, using the periodic structure of the oscillatory component
at the interface \cite{EGK2} (we must admit that such delicate
properties of oscillatory structures of solutions are known for
the 1D and radial cases only, though we expect that these
phenomena are generic), we can control the singular coefficients
in (\ref{br3}), and, in particular, to see that
\begin{equation}
 \label{f1loc}
    \ln |f| \in L^1_{\rm loc} (\ren).
\end{equation}
However, for most general geometric configurations of nonlinear
eigenfunctions $f(y)$, we do not have a proper proof of
(\ref{f1loc}) or similar estimates, so our further analysis is
still essentially formal. It is worth recalling again that our
computations below are to be understood as those dealing with the
equivalent integral equations and operators, so, in particular, we
can use the powerful facts on compactness of the resolvent
$(\BB-\l I)^{-1}$ and of the adjoint one $(\BB^*-\l I)^{-1}$ in
the corresponding weighted $L^2$-spaces. Note that, in such an
equivalent integral representation, the singular term in
(\ref{br3N}) satisfying (\ref{f1loc}) makes no principal
difficulty, so the  expansion  (\ref{br3N}) makes rather usual
sense for applying standard nonlinear operator theory.

 Thus, under natural assumptions, substituting \eqref{br3} into \eqref{sf5}, for any
$k=0,1,2,3,\cdots$\,, we find that, omitting $o(n)$ terms when
necessary,
\begin{equation*}
 \tex{
    -\nabla \cdot[(1 +n \ln |f|) \nabla \D f]
    +\frac{1-\a_k n - \mu_{1,k} n^2}{4}\, y \cdot\nabla  f
    +(\a_k+ \mu_{1,k} n) f=0\,,
    }
\end{equation*}
and, rearranging terms,
\begin{equation*}
 \tex{
    -\D^2 f -n\nabla \cdot( \ln |f| \nabla \D f)
    +\frac{1}{4}\, y \cdot\nabla  f
    -\frac{\a_k n + \mu_{1,k} n^2}{4}\, y \cdot\nabla  f
    +\a_k f+ \mu_{1,k} n f=0\,.
    }
\end{equation*}
Hence, we finally have
\begin{equation}
\label{br4}
 \tex{
    \big({\bf B} + \frac{k}{4}\,I\big) f +n \big[-\nabla \cdot( \ln |f| \nabla \D f)
    -\frac{\a_k}{4}\, y \cdot\nabla  f
    + \mu_{1,k} f \big]+o(n)   =0\,,
    }
\end{equation}
which can be written in the following form:
\begin{equation}
\label{br5}
 \tex{
    \big({\bf B} + \frac{k}{4}\,I \big) f +n \mathcal{N}_k (f)+o(n)   =0\,,
     }
\end{equation}
with the operator
\begin{equation}
\label{br6}
 \tex{
    \mathcal{N}_k (f):=-\nabla \cdot( \ln |f| \nabla \D f)
    -\frac{\a_k}{4}\, y \cdot\nabla  f
    + \mu_{1,k} f\,.
    }
\end{equation}
Subsequently, as was shown in Section \ref{S3}, we have that
\begin{equation}
\label{br7}
 \tex{
    \ker\big({\bf B} + \frac{k}{4}\, I \big)= \mathrm{span\,}\{\psi_\b\}_{|\b|=k}\quad \hbox{for any}
    \quad k=0,1,2,3,\cdots,
     }
\end{equation}
where the operator
   $ {\bf B} + \frac{k}{4}\, I$
is Fredholm of index zero and
\begin{equation*}
 \tex{
    \dim \ker\big({\bf B} + \frac{k}{4}\, I\big)=  M_k  \ge 1 \quad \hbox{for any}
    \quad k=0,1,2,3,\cdots,
     }
\end{equation*}
where $M_k$ stands for the length of the vector $\{D^\b v, \,
|\b|=k\}$, so that $M_k>1$ for  $k \ge 1$.

\ssk

\noi\underline{\sc Simple eigenvalue for $k=0$}. Since 0 is a
simple eigenvalue of ${\bf B}$ when $k=0$, i.e.,
\begin{equation*}
    \ker\,{\bf B} \oplus R[{\bf B}] = L_{\rho}^2(\ren),
\end{equation*}
the study of the case $k=0$ seems to be simpler than for other
different $k$'s because the dimension of the eigenspace is
$M_0=1$. Thus, we shall describe the behaviour of solutions for
small $n>0$ and apply the classical Lyapunov--Schmidt method to
\eqref{br5} (assuming as usual some extra necessary regularity
hypothesis to be clarified later on), in order to accomplish the
branching approach as $n\downarrow 0$, in two steps, when $k=0$
and $k$ is different from $0$.

Thus, owing to Section \ref{S3}, we have already known that $0$ is
a simple eigenvalue of ${\bf B}$, i.e., $\ker\,{\bf B}=
\mathrm{span\,}\{\psi_0\}$ is one-dimensional. Hence, denoting by
$Y_0$ the complementary invariant subspace, orthogonal to
$\psi_0^*$, we set
\begin{equation*}
    f=\psi_0+V_0,
\end{equation*}
where $V_0 \in Y_0$. According to the already well known spectral
properties of the operator ${\bf B}$, we define $P_0$ and $P_1$
such that $P_0+P_1=I$, to be the projections onto $\ker\,{\bf B}$
and $Y_0$ respectively. Finally, setting
\begin{equation}
\label{br8}
 V_0:=n \Phi_{1,0} + o(n),
\end{equation}
substituting the expression for $f$ into \eqref{br5} and passing
to the limit as $n\rightarrow 0^+$ leads to a linear inhomogeneous
equation for $\Phi_{1,0}$,
\begin{equation}
\label{br9}
    {\bf B}\Phi_{1,0}=- \mathcal{N}_0 (\psi_0),
\end{equation}
since ${\bf B} \psi_0=0$. Moreover, by  Fredholm theory, $V_0 \in
Y_0$ exists if and only if the right-hand side is orthogonal to
the one dimensional kernel of the adjoint operator ${\bf B}^*$
with $\psi_0^*=1$, because of \eqref{s16}. Hence, in the topology
of the dual space $L^2$, this requires the standard orthogonality
condition:
\begin{equation}
\label{br10}
    \big\langle \mathcal{N}_0 (\psi_0), 1\big\rangle=0.
\end{equation}
Then, \eqref{br9} has a unique solution $\Phi_{1,0} \in Y_0$
determining by \eqref{br8} a bifurcation branch for small $n>0$.
In fact, the algebraic equation \ef{br10} yields  the following
explicit expression for the coefficient $\mu_{1,0}$ of the
expansion \ef{br3} of the first eigenvalue $\a_0(n)$:
\begin{equation}
\label{br11}
    \mu_{1,0}
    :=
 \tex{
    \frac{\langle \nabla \cdot( \ln |\psi_0| \nabla \D \psi_0)
    +\frac{N}{16}\,  y \cdot\nabla  \psi_0,\psi_0^*\rangle}
    {\langle \psi_0,\psi_0^*\rangle}
 }
     =
 \tex{
    \langle \nabla \cdot( \ln |\psi_0| \nabla \D \psi_0)
    +\frac{N}{16}\,  y \cdot\nabla  \psi_0,\psi_0^*\rangle.
 }
\end{equation}

\ssk

\noi\underline{\sc Multiple eigenvalues for $k \ge 1$}. For any
$k\geq 1$, we  know that
\begin{equation*}
 \tex{
    \dim \ker\big({\bf B} + \frac{k}{4} \,I\big)=  M_k >1.
    }
\end{equation*}
  Hence, we have to use the full eigenspace expansion
\begin{equation}
\label{br12}
 \tex{
    f=\sum\limits_{|\b|=k} c_\b \hat{\psi}_\b +V_k,
 }
\end{equation}
for every $k\geq 1$. Currently, for convenience, we denote
$\{\hat{\psi}_\b\}_{|\b|=k}=\{\hat \psi_1,...,\hat \psi_{M_k}\}$
the natural basis of the $M_k$-dimensional eigenspace
$\ker\big({\bf B} + \frac{k}{4}\, I\big)$ and set
 $\psi_k =
\sum_{|\b|=k} c_\b \hat{\psi}_\b$.
 Moreover,
$V_k \in Y_k$ and $V_k=\sum_{|\b|>k} c_\b {\psi}_\b$, where $Y_k$
is the complementary invariant subspace of $\ker\big({\bf B} +
\frac{k}{4}\, I\big)$. Furthermore, in the same way, as we did for
the case $k=0$, we define the $P_{0,k}$ and $P_{1,k}$, for every
$k\geq 1$, to be the projections of $\ker\big({\bf B} +
\frac{k}{4}\, I\big)$ and $Y_k$ respectively. We also expand $V_k$
 as
\begin{equation}
\label{br13}
    V_k:=n \Phi_{1,k} + o(n).
\end{equation}
Subsequently, substituting \eqref{br12} into \eqref{br5} and
passing to the limit as $n\downarrow 0^+$, we obtain the following
equation:
\begin{equation}
\label{br14}
 \tex{
    \big({\bf B}+ \frac{k}{4}\,I\big)\Phi_{1,k}=- \mc{N}_k \big(\sum_{|\b|=k} c_\b {\psi}_\b\big),
    }
\end{equation}
under the natural ``normalizing" constraint
\begin{equation}
\label{br15}
 \tex{
    \sum\limits_{|\b|=k} c_\b=1 \quad (c_\b \ge 0).
    }
\end{equation}
Therefore, applying the Fredholm alternative, $V_k \in Y_k$ exists
if and only if the term on the right-hand side of \eqref{br14} is
orthogonal to $\ker\,\big({\bf B}+ \frac{k}{4}\,I\big)$.
Multiplying the right-hand side of \eqref{br14} by $\psi_\b^*$,
for every $|\b|=k$,  in the topology of the dual space $L^2$, we
obtain an algebraic system of $M_k+1$ equations and the same
number of unknowns, $\{c_\b, \, |\b|=k\}$ and $\mu_{1,k}$:
\be
 \label{alg1}
  \tex{
\big\langle \mc{N}_k (\sum_{|\b|=k} c_\b {\psi}_\b), \psi^*_\b
\big\rangle=0 \quad \mbox{for all} \quad  |\b|=k,
 }
 \ee
 which is indeed the Lyapunov--Schmidt branching equation
 \cite{VainbergTr}.
In general, such algebraic system are assumed to allow us to
obtain the branching parameters and hence establish the number of
different solutions induced on the given $M_k$-dimensional
eigenspace as the kernel of the operator involved.

However, we must admit and urge that the algebraic system
(\ref{alg1}) is a truly difficult issue. One of the main features
of it is as follows:
 \be
 \label{alg2}
  \fbox{$
   \mbox{(\ref{alg1}) is not variational.}
    $}
    \ee
In other words, one cannot use for (\ref{alg1}) the classic
category-genus theory of calculus of variation \cite{Berger, KZ},
to claim that the category of the kernel (equal to $M_k$) is the
least number of different critical points and hence of different
solutions.

To see (\ref{alg2}), it suffices to note that, due to (\ref{s14})
and (\ref{s16}), the generalized Hermite polynomials $\psi_\b^*$
have nothing common in the algebraic sense with the eigenfunctions
$\psi_\b$ in the $L^2$-scalar products in (\ref{alg1}).

\subsection{A digression to Hermite classic  self-adjoint theory}
 It is
worth mentioning that, for classic second-order Hermite operator
 \be
 \label{He1}
  \tex{
 \BB = \D + \frac 12\, y \cdot \n + \frac N2 \, I
 \quad \big(\mbox{then, in the $L^2$-metric,} \quad \BB^*= \D- \frac 12 \, y \cdot \n\big),
 }
 \ee
 \ef{alg2} is not the case. Indeed, by classic self-adjoint theory
 \cite[p.~48]{BS}, these eigenfunctions are related to each other
 by
   \be
   \label{He2}
  \psi_\b(y) = D^\b F(y) \equiv H_\b(y) \, F(y) \whereA F(y)=(4 \pi)^{-\frac N2}\,
  {\mathrm e}^{-\frac{|y|^2}4}
   \ee
  is the Gaussian kernel and $H_\b(y)$ are standard Hermite
  polynomials, which also define the adjoint eigenfunctions:
   \be
    \label{He3}
     \tex{
   \psi_\b^*(y)=b_\b H_\b(y) \equiv \frac{ b_\b}{F(y)} \, \psi_\b(y),
   }
     \ee
    where $b_\b$ are normalization constants.
One knows that this is a result of the symmetry of the operator
(\ref{He1}) in the weighted metric of $L^2_\rho(\ren)$, where
 $$
 \tex{
 \rho(y)= {\mathrm e}^{\frac{|y|^2}4} \sim \frac 1{F(y)}
 \LongA \BB= \frac 1\rho\,  \n \cdot (\rho \n) + \frac N2 \, I,
  \,\,\,\, \mbox{so} \,\,\,\, (\BB)_{L^2_\rho}^*=\BB.
 }
 $$

In view of the relations (\ref{He2}) and (\ref{He3}) of the
bi-orthonormal bases $\{\psi_\b\}$ and $\{\psi_\b^*\}$, the
corresponding algebraic systems such as (\ref{alg1}) can be
variational. Moreover,  even the original nonlinear elliptic
equation similar to (\ref{sf5}), where the 4th-order operator is
replaced by a natural 2nd-order one of the porous medium type:
 $$
 -\n(|f|^n \n \D f) \mapsto \n(|f|^n \n f),
  $$
  then becomes variational itself. Thus, in this case, both
  branching (local phenomena) and global extensions of
  $n$-bifurcation branches can be performed on the basis of
   powerful Lusternik--Schnirel'man category variational theory
from 1920s \cite[\S~56]{KZ}, so that existence and multiplicity
(at least, not less than in the linear case $n=0$) of solutions
are guaranteed.

\subsection{Computations for branching of dipole solutions in 2D}

To avoid excessive computations and as a self-contained example,
we now ascertain some expressions for those coefficients in the
case when $|\b|=1$, $N=2$, and $M_1=2$, so that, in our notations,
$\{\psi_\b\}_{|\b|=1}=\{\hat \psi_1, \hat \psi_2\}$. Consequently,
in this case, we obtain the following algebraic system: the
expansion coefficients of $\psi_1=c_1 \hat \psi_1+c_2 \hat \psi_2$
satisfy
\begin{equation}
\label{br16}
    \left\{\begin{array}{l}
    c_1  \langle \hat \psi_1^*,h_1 \rangle- \frac{c_1 \a_1}{4}\,
     \langle \hat \psi_1^*,y \cdot\nabla  \hat{\psi}_1 \rangle +c_1 \mu_{1,1}
    +c_2  \langle \hat \psi_1^*,h_2 \rangle- \frac{c_2 \a_1}{4}\,
     \langle \hat \psi_1^*,y \cdot\nabla  \hat{\psi}_2 \rangle = 0,\\
    c_1  \langle \hat \psi_2^*,h_1 \rangle- \frac{c_1 \a_1}{4}\,
     \langle \hat \psi_2^*,y \cdot\nabla  \hat{\psi}_1 \rangle
    + c_2  \langle \hat \psi_2^*,h_2 \rangle- \frac{c_2 \a_1}{4}\,
     \langle \hat \psi_2^*,y \cdot\nabla  \hat{\psi}_2 \rangle+ c_2 \mu_{1,1}=0,\\
    c_1+c_2=1,
    \end{array}\right.
\end{equation}
where
\begin{equation*}
   h_1:= -\nabla \cdot  [ \ln (c_1 \hat{\psi}_1+c_2\hat{\psi}_2) \nabla \D
   \hat{\psi}_1 ],\,\,
    h_2:= -\nabla \cdot
     [ \ln (c_1 \hat{\psi}_1+c_2\hat{\psi}_2) \nabla \D \hat{\psi}_2 ],
\end{equation*}
and, $c_1$, $c_2$, and $\mu_{1,1}$ are the coefficients that we
want to calculate, $\a_1$ is regarded as the value of the
parameter $\a$ denoted by \eqref{bf4} and dependent on the
eigenvalue $\l_1$, for which $\hat \psi_{1,2}$ are the associated
eigenfunctions, and $\hat \psi_{1,2}^*$ the corresponding adjoint
eigenfunctions. Hence, substituting the expression $c_2=1-c_1$
from third equation into the other two, we have the following
nonlinear algebraic system
\begin{equation}
\label{br17}
    \left\{\begin{array}{l}
    0=N_1(c_1,\mu_{1,1})- c_1\frac{\a_1}{4} \,\big[
     \langle \hat \psi_1^*,y \cdot\nabla  \hat{\psi}_1 \rangle
    -  \langle \hat \psi_1^*,y \cdot\nabla  \hat{\psi}_2 \rangle\big], \ssk\\
    0=N_2(c_1,\mu_{1,1}) - c_1\frac{\a_1}{4}\,\big[
     \langle \hat \psi_2^*,y \cdot\nabla  \hat{\psi}_1 \rangle
    - \langle \hat \psi_2^*,y \cdot\nabla  \hat{\psi}_2 \rangle\big]+ \mu_{1,1},
    \end{array}\right.
\end{equation}
where
\begin{equation*}
  \begin{split} &
   \tex{
   N_1(c_1,\mu_{1,1}):= c_1  \langle \hat \psi_1^*,h_1 \rangle
   + \langle \hat \psi_1^*,h_2 \rangle
   - \frac{\a_1}{4}\,
     \langle \hat \psi_1^*,y \cdot\nabla  \hat{\psi}_2 \rangle
   -c_1  \langle \hat \psi_1^*,h_2 \rangle +c_1 \mu_{1,1},
 }
   \\ &
    \tex{
   N_2(c_1,\mu_{1,1}):=c_1  \langle \hat \psi_2^*,h_1 \rangle
   +  \langle \hat \psi_2^*,h_2 \rangle
   - \frac{\a_1}{4} \,\langle \hat \psi_2^*,y \cdot\nabla  \hat{\psi}_2 \rangle-
   c_1 \langle \hat \psi_2^*,h_2 \rangle-c_1 \mu_{1,1}
    }
   \end{split}
\end{equation*}
represent the nonlinear parts of the algebraic system, with $h_0$
and $h_1$  depending on $c_1$.
\par

Subsequently, to guarantee existence of solutions of the system
\eqref{br16}, we apply the Brouwer Fixed Point Theorem to
\eqref{br17} by supposing that the values $c_1$ and $\mu_{1,1}$
are the unknowns, in a disc sufficiently big
$D_R(\hat{c}_1,\hat{\mu}_{1,1})$ centered in a possible
nondegenerate zero $(\hat{c}_1,\hat{\mu}_{1,1})$. Thus, we write
the system \eqref{br17} in the matrix  form
\begin{equation*}
    \binom{0}{0}= \left(\begin{array}{cc}
    -\frac{\a_1}{4}\,\big[
   \langle \hat \psi_1^*,y \cdot\nabla  \hat{\psi}_1 \rangle
  -  \langle \hat \psi_1^*,y \cdot\nabla  \hat{\psi}_2 \rangle\big] & 0\\
  -\frac{\a_1}{4} \,\big[ \langle \hat \psi_2^*,y \cdot\nabla  \hat{\psi}_1 \rangle
  - \langle \hat \psi_2^*,y \cdot\nabla  \hat{\psi}_2 \rangle\big]
    & 1\end{array}\right)\binom{c_1}{\mu_{1,1}} + \binom{N_1(c_1,\mu_{1,1})}
    {N_2(c_1,\mu_{1,1})}.
\end{equation*}
Hence, we have that the zeros of the operator
\begin{equation*}
    \mathcal{F}(c_1,\mu_{1,1}):= \mf{M} \binom{c_1}{\mu_{1,1}} + \binom{N_1(c_1,\mu_{1,1})}
    {N_2(c_1,\mu_{1,1})}
\end{equation*}
are the possible solutions of \eqref{br17}, where $\mf{M}$ is the
matrix corresponding to the linear part of the system, while
 $$
  (N_1(c_1,\mu_{1,1}),N_2(c_1,\mu_{1,1}))^T,
   $$
 corresponds to the
nonlinear part. The application $\mathcal{H}:\mathcal{A} \times
[0,1] \to \re$, defined by
\begin{equation*}
    \mathcal{H}(c_1,\mu_{1,1},t):= \mf{M}
    \binom{c_1}{\mu_{1,1}} + t\binom{N_1(c_1,\mu_{1,1})}
    {N_2(c_1,\mu_{1,1})},
\end{equation*}
provides us with a homotopy transformation from the function
$\mathcal{F}(c_1,\mu_{1,1})= \mathcal{H}(c_1,\mu_{1,1},1)$ to its
linearization
\begin{equation}
\label{br18}
    \mathcal{H}(c_1,\mu_{1,1},0):= \mf{M} \binom{c_1}{\mu_{1,1}}.
\end{equation}

Thus, the system \eqref{br17} possesses a nontrivial solution if
\eqref{br18} has a nondegenerate zero, in other words, if the next condition
is satisfied
\begin{equation}
\label{br19}
     \langle \hat \psi_1^*,y \cdot\nabla  \hat{\psi}_1 \rangle-
     \langle \hat \psi_1^*,y \cdot\nabla  \hat{\psi}_2 \rangle\neq 0.
\end{equation}
Note that, if the substitution would have been $c_1=1-c_2$, the
condition might also be
\begin{equation*}
     \langle \hat \psi_2^*,y \cdot\nabla  \hat{\psi}_2 \rangle-
     \langle \hat \psi_2^*,y \cdot\nabla  \hat{\psi}_1 \rangle\neq 0.
\end{equation*}
Then, under condition \eqref{br19}, the system \eqref{br17} can be
written in the form
\begin{equation}
\label{br20}
    \binom{c_1-\hat{c}_1}{\mu_{1,1}-\hat{\mu}_{1,1}}=-\mathcal{M}^{-1}
    \binom{N_1(c_1,\mu_{1,1})-\hat{c}_1}
    {N_2(c_1,\mu_{1,1})- \hat{\mu}_{1,1}},
\end{equation}
which can be interpreted as a fixed point equation. Moreover,
applying Brouwer's Fixed Point Theorem, we have that
\begin{align*}
    \hbox{Ind}((\hat{c}_1,\hat{\mu}_{1,1}),\mathcal{H}(\cdot,\cdot,0)) & =
    \mathcal{Q}_{C_R (\hat{c}_1,\hat{\mu}_{1,1})}(\mathcal{H}(\cdot,\cdot,0))\\ &
    =\hbox{Deg}(\mathcal{H}(\cdot,\cdot,0),D_R (\hat{c}_1,\hat{\mu}_{1,1}))\\ &
    = \hbox{Deg}(\mathcal{F}(c_1,\mu_{1,1}), D_R (\hat{c}_1,\hat{\mu}_{1,1})),
\end{align*}
where $\mathcal{Q}_{C_R
(\hat{c}_1,\hat{\mu}_{1,1})}(\mathcal{H}(\cdot,\cdot,0))$ defines
the number of rotations of the function
$\mathcal{H}(\cdot,\cdot,0)$ around the curve $C_R
(\hat{c}_1,\hat{\mu}_{1,1})$ and
$\hbox{Deg}(\mathcal{H}(\cdot,\cdot,0),D_R
(\hat{c}_1,\hat{\mu}_{1,1}))$ denotes the topological degree of
$\mathcal{H}(\cdot,\cdot,0)$ in $D_R (\hat{c}_1,\hat{\mu}_{1,1})$.
Owing to classical topological methods, both are equal.
\par

Thus, once we have proved the existence of solutions, we achieve
some expressions for the coefficients required:
\begin{equation*}
    \left\{
    \begin{matrix}
    \mu_{1,1} =c_2 ( \langle \hat \psi_1^*+\hat \psi_2^*,h_1-h_2 \rangle
    - \frac{\a_1}{4}\,
     \langle \hat \psi_1^*+\psi_2^*
    ,y \cdot\nabla  \hat{\psi}_1-y \cdot\nabla  \hat{\psi}_2 \rangle)
    \ssk\\
     -  \langle \hat \psi_1^*+\hat \psi_2^*,h_1 \rangle
    + \frac{\a_1}{4}\,
     \langle \hat \psi_1^*+\hat \psi_2^*
    ,y \cdot\nabla  \hat{\psi}_1  \rangle, \qquad\qquad\qquad\,\,
    \ssk \\
    c_1
    =1-c_2.\,\,\,\qquad\qquad\qquad\qquad\qquad\qquad\qquad\qquad\qquad\qquad\qquad
      \end{matrix}
      \right.
\end{equation*}
The expressions for the coefficients in a general case might be
accomplished after some tedious calculations, otherwise similar to
those performed above. Note that, in general, those nonlinear
finite-dimensional algebraic problems are rather complicated, and
the problem of an optimal estimate of the number of different
solutions remains open. Moreover, reliable multiplicity results
are very difficult to obtain. We expect that this number should be
somehow related (and even sometimes coincides) with the dimension
of the corresponding eigenspace of the linear operators ${\bf
B}+\frac{k}{4}\, I$, for any $k=0,1,2,\ldots\,$. This is a
conjecture only that may be  too illusive; see further supportive
analysis presented below.
\par

 However, we devote the remaining of this section to a possible answer
to that conjecture, which is not totally complete though, since we
are imposing some conditions.
\par

Thus, in order to detect the number of solutions of the nonlinear
algebraic system \eqref{br16}, we proceed to reduce this system to
a single equation for one of the unknowns. As a first step,
integrating by parts in the terms in which $h_1$ and $h_2$ are
involved and rearranging terms in the first two equations of the
system \eqref{br16}, we arrive at
 $$
 \left\{
\begin{matrix}
 \tex{
   \int\limits_{\ren} \nabla \psi_1^* \cdot \ln (c_1 \hat{\psi}_1+c_2\hat{\psi}_2) \nabla \D
   (c_1 \hat{\psi}_1
 }
    \tex{
    + c_2 \hat{\psi}_2)
    } \qquad\qquad\quad
    \\
 \tex{
 - c_1 \frac{\a_1}{4}\,
    \int\limits_{\ren} \hat \psi_1^* y \cdot\nabla  \hat{\psi}_1
    +c_1 \mu_{1,1}- c_2  \frac{\a_1}{4}\,
    \int\limits_{\ren} \hat\psi_1^* y \cdot\nabla  \hat{\psi}_2 = 0,
    } \ssk\ssk
    \\
 \tex{
   \int\limits_{\ren} \nabla \hat \psi_2^*\cdot \ln (c_1 \hat{\psi}_1+c_2\hat{\psi}_2) \nabla \D
   (c_2 \hat{\psi}_1+ c_2 \hat{\psi}_2)
 } \qquad\qquad\quad
   \\
 \tex{
     -c_1  \frac{\a_1}{4}\,
   \int\limits_{\ren} \hat \psi_2^* y \cdot\nabla  \hat{\psi}_1
    + c_2 \mu_{1,1}
 }
      \tex{
    -c_2  \frac{\a_1}{4}\,
    \int\limits_{\ren} \hat\psi_2^* y \cdot\nabla  \hat{\psi}_2 =0.
     }
\end{matrix}
 \right.
 $$
 By the third equation, we have that $c_1 = 1- c_2$, and   hence,
setting $c_1 \hat{\psi}_1+c_2\hat{\psi}_2 =
\hat{\psi}_1+(\hat{\psi}_2- \hat{\psi}_1)c_2$ and substituting
these into those new expressions for the first two equations of
the system, we find that
 \be
 \label{br59}
 \left\{
\begin{matrix}
 \tex{
   \int\limits_{\ren} \nabla \hat \psi_1^* \cdot
   \ln (\hat{\psi}_1+
 }
     \tex{
     (\hat{\psi}_2-\hat{\psi}_1)c_2) \nabla \D
   (\hat{\psi}_1+(\hat{\psi}_2- \hat{\psi}_1)c_2) +\mu_{1,1}
    -c_2 \mu_{1,1}
     }
    \\ 
 \tex{
     - \frac{\a_1}{4}\,
    \int\limits_{\ren} \hat \psi_1^* y \cdot\nabla  \hat{\psi}_1
    + c_2  \frac{\a_1}{4}\, \int\limits_{\ren} \hat{\psi}_1^* y
    \cdot(\nabla  \hat{\psi}_1-\nabla \hat{\psi}_2)  = 0,
     } \qquad\qquad\quad\,\, \ssk\ssk
      \\
 \tex{
   \int\limits_{\ren} \nabla \hat \psi_2^*\cdot
   \ln (\hat{\psi}_1+
   }
     \tex{
     (\hat{\psi}_2-\hat{\psi}_1)c_2) \nabla \D
   (\hat{\psi}_1+(\hat{\psi}_2- \hat{\psi}_1)c_2)
    + c_2 \mu_{1,1}
 }\qquad\,\,\,
     \\  
 \tex{
     - \frac{\a_1}{4}\,
    \int\limits_{\ren} \hat \psi_2^* y \cdot\nabla  \hat{\psi}_1
    +c_2  \frac{\a_1}{4}\,
    \int\limits_{\ren} \hat \psi_2^* y
    \cdot(\nabla  \hat{\psi}_1-\nabla \hat{\psi}_2) =0.
    }\qquad\qquad\quad\,\,
\end{matrix}
 \right.
\end{equation}
Subsequently, adding both equations, we have that
\begin{align*}
   \mu_{1,1} &  =
 \tex{
    - \int\limits_{\ren}  (\nabla \hat \psi_1^*+ \nabla \hat \psi_2^*) \cdot
   \ln (\hat{\psi}_1+  (\hat{\psi}_2-\hat{\psi}_1)c_2) \nabla \D
   (\hat{\psi}_1+(\hat{\psi}_2- \hat{\psi}_1)c_2)
   }
    \\ &
 \tex{
     + \frac{\a_1}{4}\,
    \int\limits_{\ren} (\psi_1^*+ \psi_2^*) y
    \cdot \nabla  \hat{\psi}_1
    - c_2  \frac{\a_1}{4}\,
    \int\limits_{\ren} (\hat \psi_1^*+ \hat \psi_2^*) y
    \cdot(\nabla  \hat{\psi}_2-\nabla \hat{\psi}_1).
    }
\end{align*}
Thus, substituting it into the second equation of \eqref{br59}, we
obtain the following equation with the single unknown $c_2$:
\begin{equation}
\label{br60}
\begin{split}
   &
    \tex{
    -c_2^2 \frac{\a_1}{4}\,
    \int\limits_{\ren} (\hat \psi_1^*+ \hat \psi_2^*) y
    \cdot (\nabla  \hat{\psi}_2-\nabla \hat{\psi}_1) +
    c_2 \frac{\a_1}{4}(\,
    \int\limits_{\ren} (\hat\psi_1^*+ 2 \hat\psi_2^*) y
    \cdot \nabla  \hat{\psi}_1- \,
    \int\limits_{\ren} \hat\psi_2^* y \cdot \nabla \hat{\psi}_2)
    }
    \\ &
 \tex{
     - \frac{\a_1}{4}\,
    \int\limits_{\ren} \hat\psi_2^* y \cdot\nabla  \hat{\psi}_1
    +\int\limits_{\ren} \nabla \psi_2^*\cdot
   \ln (\hat{\psi}_1+ (\hat{\psi}_2-\hat{\psi}_1)c_2) \nabla \D
   (\hat{\psi}_1+(\hat{\psi}_2- \hat{\psi}_1)c_2)
   }
   \\ &
 \tex{
    -c_2 \int\limits_{\ren} (\nabla \hat \psi_1^* +\nabla \hat\psi_2^*)\cdot
   \ln (\hat{\psi}_1+ (\hat{\psi}_2-\hat{\psi}_1)c_2) \nabla \D
   (\hat{\psi}_1+(\hat{\psi}_2- \hat{\psi}_1)c_2) =0,
   }
\end{split}
\end{equation}
which can be written in the following way:
\begin{equation}
 \label{FF1}
    c_2^2 A + c_2 B + C + \o(c_2) \equiv \mf{F} (c_2) + \o
    (c_2)=0.
\end{equation}
 Here, $\o(c_2)$ can be considered as perturbation of the quadratic
form $\mf{F} (c_2)$ with the coefficients
defined by
\begin{align*}
    &
     \tex{
     A:=  -\frac{\a_1}{4}\,
    \int\limits_{\ren} (\hat\psi_1^*+ \hat\psi_2^*) y
    \cdot (\nabla  \hat{\psi}_2-\nabla \hat{\psi}_1),
 }
    \\ &
    \tex{
  B:= \frac{\a_1}{4}(\,
    \int\limits_{\ren} (\hat\psi_1^*+ 2 \hat\psi_2^*) y
    \cdot \nabla  \hat{\psi}_1- \,
    \int\limits_{\ren} \hat\psi_2^* y \cdot \nabla \hat{\psi}_2),
    }
    , \quad
  C:=
 \tex{
   - \frac{\a_1}{4}\,
    \int\limits_{\ren} \hat\psi_2^* y \cdot\nabla  \hat{\psi}_1,
    }
     \\ &
  \o(c_2):=
  \tex{
   \int\limits_{\ren} \nabla \hat\psi_2^*\cdot
   \ln (\hat{\psi}_1+ (\hat{\psi}_2-\hat{\psi}_1)c_2) \nabla \D
   (\hat{\psi}_1+(\hat{\psi}_2- \hat{\psi}_1)c_2)
   }
   \\ &
    \tex{
    -c_2 \int\limits_{\ren} (\nabla \hat\psi_1^*+\nabla \hat\psi_2^*)\cdot
   \ln (\hat{\psi}_1+ (\hat{\psi}_2-\hat{\psi}_1)c_2) \nabla \D
   (\hat{\psi}_1+(\hat{\psi}_2- \hat{\psi}_1)c_2).
   }
\end{align*}

Since, due to the normalizing constraint \eqref{br15}, $c_2 \in
[0,1]$, solving the quadratic equation $ \mf{F} (c_2)$ yields:
\begin{enumerate}
\item[(i)] $c_2=0 \Longrightarrow \mf{F}(0) = C $;
\item[(ii)] $c_2 =1 \Longrightarrow \mf{F}(1)= A+B+C $; and
\item[(iii)] differentiating $\mf{F}$ with respect to $c_2$, we obtain that
$\mf{F}'(c_2) = 2 c_2 A+B$. Then, the critical point of the
function $\mf{F}$ is $c_2^* = -\frac{B}{2A}$ and its image is $
\mf{F}(c_2^*)=- \frac{B}{4A}+C$.
\end{enumerate}

Consequently, the conditions that must be imposed in order to have
more than one solution (we already know the existence of at least
one solution) are as follows:
\begin{enumerate}
\item[(a)] $C  (A+B+C) >0$;
\item[(b)] $C \big(-\frac{B}{4A}+C\big)<0 $; and
\item[(c)] $0<-\frac{B}{2A}<1$.
\end{enumerate}
Note that, for $-\frac{B}{4A}+C= 0$, we have just a single
solution.
\par
Hence, considering the equation again in the form
 $$
 \mf{F} (c_2) + \o
(c_2)=0,
 $$
  where $\o (c_2)$ is a perturbation of the quadratic form
$\mf{F} (c_2)$, and bearing in mind that the objective is to
detect  the number of solutions of the system \eqref{br16}, we
need to control  somehow this perturbation. Under the conditions
(a), (b), and (c), $\mf{F} (c_2)$ possesses exactly two solutions.
Therefore, controlling the possible oscillations of the
perturbation $\o (c_2)$ in such a way that
\begin{equation*}
    \left\| \o (c_2) \right\|_{L^\infty} \leq \mf{F}(c_2^*),
\end{equation*}
we can assure that the number of solutions for \eqref{br16} is
exactly two. This is  the dimension of the kernel of the operator
${\bf B}+\frac{1}{4}\, I$ (as we expected in our more general
conjecture).

 The above particular example shows how
difficult are the questions on existence and multiplicity of
solutions for such non-variational branching problems.
 Recall that the actual values of the coefficients $A$, $B$, $C$,
 and others, which the number of solutions crucially depend on, is
 difficult even estimate numerically in view of a complicated
 nature of the eigenfunctions \ef{s14} involved, to say nothing of
 the nonlinear perturbation $\o(c_2)$.

\subsection{Branching computations for $|\b|=2$}

\ssk

\ssk

Overall, the above analysis provides us with some expressions for
the solutions for the self-similar equation \eqref{sf5} depending
on the value of $k$. Actually, we can achieve those expressions
for every critical value $\a_k$, but again the calculus gets
rather difficult. For the sake of completeness, we now analyze the
case $|\b|=2$ and $M_2=3$, so that $\{\psi_\b\}_{|\b|=2}=\{\hat
\psi_1, \hat \psi_2, \hat \psi_3\}$ stands for a basis of the
eigenspace $\ker\big({\bf B} + \frac{1}{2}\, I \big)$, with $k=2$
($\l_k=- \frac k4$).

Thus, in this case, performing in a similar way as was done for
\eqref{br16} with $\psi_2= c_1 \hat \psi_1+c_2 \hat \psi_2+c_3
\hat \psi_3$, we arrive at the following algebraic system:
\begin{equation}
\label{br61}
    \left\{\begin{array}{l}
    \begin{array}{r}
     c_1  \langle \hat \psi_1^*,h_1 \rangle
    +c_2  \langle \hat \psi_1^*,h_2 \rangle
    +c_3  \langle \hat \psi_1^*,h_3 \rangle - \frac{c_1 \a_2}{4}\,
     \langle \hat \psi_1^*,y \cdot\nabla  \hat{\psi}_1 \rangle
    - \frac{c_2 \a_2}{4}\,
     \langle \hat \psi_1^*,y \cdot\nabla  \hat{\psi}_2 \rangle
    \\- \frac{c_3 \a_2}{4}\,
     \langle \hat \psi_1^*,y \cdot\nabla  \hat{\psi}_3 \rangle +c_1 \mu_{1,2}= 0,
    \ssk \\
     c_1  \langle \hat \psi_2^*,h_1 \rangle
    + c_2  \langle \hat \psi_2^*,h_2 \rangle
    + c_2  \langle \hat \psi_2^*,h_3 \rangle  - \frac{c_1
    \a_2}{4}\,
     \langle \hat \psi_2^*,y \cdot\nabla  \hat{\psi}_1 \rangle
    - \frac{c_2 \a_2}{4}\,
     \langle \hat \psi_2^*,y \cdot\nabla  \hat{\psi}_2 \rangle
     \\  - \frac{c_3 \a_2}{4}\,
     \langle \hat \psi_2^*,y \cdot\nabla  \hat{\psi}_3 \rangle+ c_2 \mu_{1,2}=0,
     \ssk \\
     c_1  \langle \hat \psi_3^*,h_1 \rangle
    + c_2  \langle \hat \psi_3^*,h_2 \rangle
    + c_2  \langle \hat \psi_3^*,h_3 \rangle - \frac{c_1
    \a_2}{4}\,
     \langle \hat \psi_3^*,y \cdot\nabla  \hat{\psi}_1 \rangle
    - \frac{c_2 \a_2}{4}\,
     \langle \hat \psi_3^*,y \cdot\nabla  \hat{\psi}_2 \rangle
     \\  - \frac{c_3 \a_2}{4}\,
     \langle \hat \psi_3^*,y \cdot\nabla  \hat{\psi}_3 \rangle+ c_3 \mu_{1,2}=0,
     \end{array}\\
    c_1+c_2+c_3=1,
    \end{array}\right.
\end{equation}
where
 $$
  \begin{matrix}
   h_1:= -\nabla \cdot  [ \ln (c_1 \hat{\psi}_1+c_2\hat{\psi}_2
   +c_3 \hat{\psi}_3) \nabla \D
   \hat{\psi}_1 ], \\
    h_2:= -\nabla \cdot
     [ \ln (c_1 \hat{\psi}_1+c_2\hat{\psi}_2
     +c_3 \hat{\psi}_3) \nabla \D \hat{\psi}_2 ], \\
   h_3:= -\nabla \cdot  [ \ln (c_1 \hat{\psi}_1+c_2\hat{\psi}_2
   +c_3 \hat{\psi}_3) \nabla \D
   \hat{\psi}_3 ],
\end{matrix}
 $$
and $c_1$, $c_2$, $c_3$, and $\mu_{1,2}$ are the unknowns to be
evaluated. Also, $\a_2$ is regarded as the value of the parameter
$\a$ denoted by \eqref{bf4} and is dependent on the eigenvalue
$\l_2$ with $\hat{\psi}_1,\hat{\psi}_2,\hat{\psi}_3$ representing
the associated eigenfunctions and
$\hat{\psi}_1^*,\hat{\psi}_2^*,\hat{\psi}_3^*$ the corresponding
adjoint eigenfunctions.
\par
Subsequently, substituting $c_3=1-c_1-c_2$ into the first three
equations and performing an argument based upon the Brouwer Fixed
Point Theorem and the topological degree as the one done above for
the case $|\b|=1$, we ascertain the existence of a nondegenerate
solution of the algebraic system if the following condition is
satisfied:
\begin{equation}
\label{br62}
     \langle \hat \psi_1^*,y \cdot\nabla  (\hat{\psi}_3-\hat{\psi}_1) \rangle
     \langle \hat \psi_2^*,y \cdot\nabla  (\hat{\psi}_3-\hat{\psi}_2) \rangle-
     \langle \hat \psi_1^*,y \cdot\nabla  (\hat{\psi}_3-\hat{\psi}_2) \rangle
     \langle \hat \psi_2^*,y \cdot\nabla  (\hat{\psi}_3-\hat{\psi}_1) \rangle\neq 0.
\end{equation}
Note that, by similar substitutions, other conditions might be
obtained.
\par
Furthermore, once we know the existence of at least one solution,
we proceed now with a possible way of computing the number of
solutions of the nonlinear algebraic system \eqref{br61}.
Obviously, since the dimension of the eigenspace is bigger than
that in the case when $|\b|=1$, the difficulty to obtain
multiplicity results increases.
\par
Firstly, integrating by parts in the nonlinear terms, in which
$h_1$, $h_2$ and $h_3$ are involved, and rearranging terms in the
first three equations gives
\begin{align*}
 \tex{
   \int\limits_{\ren} \nabla \psi_1^* \cdot
   \ln (c_1 \hat{\psi}_1+c_2\hat{\psi}_2+ c_3 \hat{\psi}_3)}
   & \tex{ \nabla \D
   (c_1 \hat{\psi}_1+ c_2 \hat{\psi}_2+ c_3 \hat{\psi}_3)
 }
    \tex{
     - c_1 \frac{\a_2}{4}\,
    \int\limits_{\ren} \hat \psi_1^* y \cdot\nabla  \hat{\psi}_1+c_1 \mu_{1,2}
    }
    \\ &
 \tex{
    - c_2  \frac{\a_2}{4}\,
    \int\limits_{\ren} \hat\psi_1^* y \cdot\nabla  \hat{\psi}_2
    - c_3  \frac{\a_2}{4}\,
    \int\limits_{\ren} \hat\psi_1^* y \cdot\nabla  \hat{\psi}_3 =
    0,
    }
\end{align*}
\begin{align*}
 \tex{
   \int\limits_{\ren} \nabla \hat \psi_2^*\cdot
   \ln (c_1 \hat{\psi}_1+c_2\hat{\psi}_2+ c_3 \hat{\psi}_3)} &
   \tex{ \nabla \D
   (c_2 \hat{\psi}_1+ c_2 \hat{\psi}_2+ c_3 \hat{\psi}_3)
 }
 \tex{
     -c_1  \frac{\a_2}{4}\,
   \int\limits_{\ren} \hat \psi_2^* y \cdot\nabla  \hat{\psi}_1
    + c_2 \mu_{1,2}
 }
     \\ &
      \tex{
    -c_2  \frac{\a_2}{4}\,
    \int\limits_{\ren} \hat\psi_2^* y \cdot\nabla  \hat{\psi}_2
    -c_3  \frac{\a_2}{4}\,
   \int\limits_{\ren} \hat \psi_2^* y \cdot\nabla  \hat{\psi}_3=0,
     }
\end{align*}
\begin{align*}
 \tex{
   \int\limits_{\ren} \nabla \hat \psi_3^*\cdot
   \ln (c_1 \hat{\psi}_1+c_2\hat{\psi}_2+ c_3 \hat{\psi}_3)} &
   \tex{ \nabla \D
   (c_2 \hat{\psi}_1+ c_2 \hat{\psi}_2+ c_3 \hat{\psi}_3)
 }
 \tex{
     -c_1  \frac{\a_2}{4}\,
   \int\limits_{\ren} \hat \psi_3^* y \cdot\nabla  \hat{\psi}_1
    + c_3 \mu_{1,2}
 }
     \\ &
      \tex{
    -c_2  \frac{\a_2}{4}\,
    \int\limits_{\ren} \hat\psi_3^* y \cdot\nabla  \hat{\psi}_2
    -c_3  \frac{\a_2}{4}\,
   \int\limits_{\ren} \hat \psi_3^* y \cdot\nabla  \hat{\psi}_3=0.
     }
\end{align*}
By the fourth equation, we have that $c_1 = 1-c_2-c_3$. Then,
setting
\begin{equation*}
    c_1 \hat{\psi}_1+c_2\hat{\psi}_2+ c_3 \hat{\psi}_3 =
    \hat{\psi}_1+c_2(\hat{\psi}_2-\hat{\psi}_1)+ c_3
    (\hat{\psi}_3-\hat{\psi}_1)
\end{equation*}
and substituting it into the expressions obtained above for the
first three equations of the system yield
\begin{align*}
 \tex{
   \int\limits_{\ren} \nabla \hat \psi_1^* \cdot
   \ln (\hat{\psi}_1+
 }
    &
     \tex{
     (\hat{\psi}_2-\hat{\psi}_1)c_2+(\hat{\psi}_3-\hat{\psi}_1)c_3) \nabla \D
   (\hat{\psi}_1+(\hat{\psi}_2- \hat{\psi}_1)c_2+
   (\hat{\psi}_3- \hat{\psi}_1)c_3)
     }
     \\ &
 \tex{
     +\mu_{1,2}
    -c_2 \mu_{1,2}- c_3 \mu_{1,2} - \frac{\a_2}{4}\,
    \int\limits_{\ren} \hat \psi_1^* y \cdot\nabla  \hat{\psi}_1
     }
    \\ &
 \tex{
    + \frac{\a_2}{4} \,\int\limits_{\ren} \hat{\psi}_1^* y
    \cdot( (\nabla  \hat{\psi}_1-\nabla \hat{\psi}_2)c_2 +
    (\nabla \hat{\psi}_1 - \nabla\hat{\psi}_3)c_3)  = 0,
     }
\end{align*}
\begin{equation}
\begin{split}
\label{br63}
 \tex{
   \int\limits_{\ren} \nabla \hat \psi_2^* \cdot
   \ln (\hat{\psi}_1+
 }
    &
     \tex{
     (\hat{\psi}_2-\hat{\psi}_1)c_2+(\hat{\psi}_3-\hat{\psi}_1)c_3) \nabla \D
   (\hat{\psi}_1+(\hat{\psi}_2- \hat{\psi}_1)c_2+
   (\hat{\psi}_3- \hat{\psi}_1)c_3)
     }
     \\ &
 \tex{
     +c_2 \mu_{1,2}
    - \frac{\a_2}{4}\,
    \int\limits_{\ren} \hat \psi_2^* y \cdot\nabla  \hat{\psi}_1
     }
    \\ &
 \tex{
    + \frac{\a_2}{4}\, \int\limits_{\ren} \hat{\psi}_2^* y
    \cdot( (\nabla  \hat{\psi}_1-\nabla \hat{\psi}_2)c_2 +
    (\nabla \hat{\psi}_1 - \nabla\hat{\psi}_3)c_3)  = 0,
     }
\end{split}
\ee
\begin{align*}
 \tex{
   \int\limits_{\ren} \nabla \hat \psi_3^* \cdot
   \ln (\hat{\psi}_1+
 }
    &
     \tex{
     (\hat{\psi}_2-\hat{\psi}_1)c_2+(\hat{\psi}_3-\hat{\psi}_1)c_3) \nabla \D
   (\hat{\psi}_1+(\hat{\psi}_2- \hat{\psi}_1)c_2+
   (\hat{\psi}_3- \hat{\psi}_1)c_3)
     }
     \\ &
 \tex{
     +c_3\mu_{1,2}
    - \frac{\a_2}{4}\,
    \int\limits_{\ren} \hat \psi_3^* y \cdot\nabla  \hat{\psi}_1
     }
    \\ &
 \tex{
    + \frac{\a_2}{4} \,\int\limits_{\ren} \hat{\psi}_3^* y
    \cdot( (\nabla  \hat{\psi}_1-\nabla \hat{\psi}_2)c_2 +
    (\nabla \hat{\psi}_1 - \nabla\hat{\psi}_3)c_3)  = 0.
     }
\end{align*}
Now, adding the first equation of \eqref{br63} to the other two,
we have that
\begin{align*}
 \tex{
   \int\limits_{\ren} (\nabla \hat \psi_1^*+ \nabla \hat{\psi}_2^*) \cdot
   \ln (\hat{\psi}_1+
 }
    &
     \tex{
     (\hat{\psi}_2-\hat{\psi}_1)c_2+(\hat{\psi}_3-\hat{\psi}_1)c_3) \nabla \D
   (\hat{\psi}_1+(\hat{\psi}_2- \hat{\psi}_1)c_2+
   (\hat{\psi}_3- \hat{\psi}_1)c_3)
     }
     \\ &
 \tex{
     +\mu_{1,2}
    - c_3 \mu_{1,2} - \frac{\a_2}{4}\,
    \int\limits_{\ren} (\hat \psi_1^*+ \hat{\psi}_2^*) y \cdot\nabla  \hat{\psi}_1
     }
    \\ &
 \tex{
    + \frac{\a_2}{4}\, \int\limits_{\ren} (\hat \psi_1^*+ \hat{\psi}_2^*) y
    \cdot( (\nabla  \hat{\psi}_1-\nabla \hat{\psi}_2)c_2 +
    (\nabla \hat{\psi}_1 - \nabla \hat{\psi}_3)c_3)  = 0,
     }
\end{align*}
\begin{align*}
 \tex{
   \int\limits_{\ren} (\nabla \hat \psi_1^*+ \nabla \hat{\psi}_3^*) \cdot
   \ln (\hat{\psi}_1+
 }
    &
     \tex{
     (\hat{\psi}_2-\hat{\psi}_1)c_2+(\hat{\psi}_3-\hat{\psi}_1)c_3) \nabla \D
   (\hat{\psi}_1+(\hat{\psi}_2- \hat{\psi}_1)c_2+
   (\hat{\psi}_3- \hat{\psi}_1)c_3)
     }
     \\ &
 \tex{
     +\mu_{1,2}
    - c_2 \mu_{1,2} - \frac{\a_2}{4}\,
    \int\limits_{\ren} (\hat \psi_1^*+ \hat{\psi}_3^*) y \cdot\nabla  \hat{\psi}_1
     }
    \\ &
 \tex{
    + \frac{\a_2}{4}\, \int\limits_{\ren} (\hat \psi_1^*+ \hat{\psi}_3^*) y
    \cdot( (\nabla  \hat{\psi}_1-\nabla \hat{\psi}_2)c_2 +
    (\nabla \hat{\psi}_1 - \nabla \hat{\psi}_3)c_3)  = 0.
     }
\end{align*}
Subsequently, subtracting those equations yields
\begin{align*}
 \tex{
    \mu_{1,2}
 }
 &
 \tex{
     = \frac{1}{c_3-c_2} \,\big[ \int\limits_{\ren}
      (\nabla \hat \psi_2^*- \nabla \hat{\psi}_3^*) \cdot
   \ln \Psi \nabla \D \Psi -\frac{\a_2}{4}\,
    \int\limits_{\ren} (\hat \psi_2^*- \hat{\psi}_3^*) y \cdot\nabla  \hat{\psi}_1
  }
     \\ &
  \tex{
    +\frac{\a_2}{4}\, \int\limits_{\ren} (\hat \psi_2^*- \hat{\psi}_3^*) y
    \cdot( (\nabla  \hat{\psi}_1-\nabla \hat{\psi}_2)c_2 +
    (\nabla \hat{\psi}_1 - \nabla \hat{\psi}_3)c_3) \big],
  }
\end{align*}
where $\Psi = \hat{\psi}_1+
(\hat{\psi}_2-\hat{\psi}_1)c_2+(\hat{\psi}_3-\hat{\psi}_1)c_3$.
Thus, substituting it into \eqref{br63} (note that, from the
substitution into one of the last two equations, we obtain the
same equation), we arrive at the following system, with $c_2$ and
$c_3$ as the unknowns:
\begin{align*}
 \tex{
   c_3 \int\limits_{\ren} (\nabla \hat{\psi}_1^*-
   \nabla \hat{\psi}_2^*  }
    &
     \tex{
     +\nabla \hat{\psi}_3^*) \cdot
   \ln \Psi
     \nabla \D \Psi - c_2 \int\limits_{\ren} (\nabla \hat{\psi}_1^*+
   \nabla \hat{\psi}_2^*-\nabla \hat{\psi}_3^*) \cdot
   \ln \Psi \nabla \D \Psi} \\ & \tex{ + \int\limits_{\ren} (
   \nabla \hat{\psi}_2^*-\nabla \hat{\psi}_3^*) \cdot
   \ln \Psi- \frac{\a_2}{4} \int\limits_{\ren} (\hat{\psi}_2^*-\hat{\psi}_3^*)
   y \cdot\nabla \hat{\psi}_1 \,} \\ & \tex{
    +c_2 \frac{\a_2}{4}\,
    [\int\limits_{\ren} (\hat{\psi}_2^*-\hat{\psi}_3^*) y \cdot\nabla
    (2\hat{\psi}_1 -\hat{\psi}_2)- \int\limits_{\ren}
    \hat{\psi}_1^* y \cdot\nabla  \hat{\psi}_1]
     }
     \\ &
 \tex{
    + c_3 \frac{\a_2}{4}\, [\int\limits_{\ren} (\hat{\psi}_2^*-\hat{\psi}_3^*) y \cdot\nabla
    (2\hat{\psi}_1 -\hat{\psi}_3)- \int\limits_{\ren}
    \hat{\psi}_1^* y \cdot\nabla  \hat{\psi}_1]
     }
      \end{align*}
      \begin{align*}
  & \tex{
     + c_2 c_3  \frac{\a_2}{4}[\, \int\limits_{\ren} \hat{\psi}_1^* y
    \cdot( \nabla  \hat{\psi}_3-\nabla \hat{\psi}_2) - \, \int\limits_{\ren}
    (\nabla \hat{\psi}_2^* - \nabla\hat{\psi}_3^*) y \cdot (2\nabla \hat{\psi}_1 -
    \nabla \hat{\psi}_2-\nabla \hat{\psi}_3)]} \\ &
    \tex{
    + c_3^2  \frac{\a_2}{4}\, \int\limits_{\ren} (\hat{\psi}_1^*-
    \hat{\psi}_2^*+\hat{\psi}_3^*) y
    \cdot( \nabla \hat{\psi}_1 - \nabla\hat{\psi}_3)}\\ &
    \tex{
    - c_2^2  \frac{\a_2}{4}\, \int\limits_{\ren} (\hat{\psi}_1^*+
    \hat{\psi}_2^*-\hat{\psi}_3^*) y
    \cdot(\nabla  \hat{\psi}_1-\nabla \hat{\psi}_2)     = 0,
     }
\end{align*}
\begin{align*}
 \tex{
   c_3 \int\limits_{\ren} \nabla \hat{\psi}_2^* \cdot
   \ln \Psi
 }
    &
     \tex{
     \nabla \D \Psi - c_2 \int\limits_{\ren} \nabla \hat{\psi}_3^* \cdot
   \ln \Psi \nabla \D \Psi - c_3 \frac{\a_2}{4}\,
    \int\limits_{\ren} \hat{\psi}_2^* y \cdot\nabla  \hat{\psi}_1
    +c_2 \frac{\a_2}{4}\,
    \int\limits_{\ren} \hat{\psi}_3^* y \cdot\nabla  \hat{\psi}_1
     }
     \\ &
 \tex{
    + c_3 \frac{\a_2}{4}\, \int\limits_{\ren} \hat{\psi}_2^* y
    \cdot( (\nabla  \hat{\psi}_1-\nabla \hat{\psi}_2)c_2 +
    (\nabla \hat{\psi}_1 - \nabla\hat{\psi}_3)c_3)
     }
     \\ &
 \tex{
    - c_2 \frac{\a_2}{4}\, \int\limits_{\ren} \hat{\psi}_3^* y
    \cdot( (\nabla  \hat{\psi}_1-\nabla \hat{\psi}_2)c_2 +
    (\nabla \hat{\psi}_1 - \nabla\hat{\psi}_3)c_3)  = 0.
     }
\end{align*}
These can be re-written in the following form:
\begin{equation}
\label{br64}
    \begin{split}
    &
    A_1 c_2^2+ B_1 c_3^2 + C_1 c_2 +D_1 c_3+ E_1 c_2 c_3 +\o_1 (c_2,c_3)=0,
    \\ &
    A_2 c_2^2+ B_2 c_3^2 + C_2 c_2 +D_2 c_3+ E_2 c_2 c_3 +\o_2 (c_2,c_3)=0,
    \end{split}
\end{equation}
 where
\begin{align*}
    \tex{
    \o_1 (c_2,c_3)} & \tex{
   := c_3 \int\limits_{\ren} (\nabla \hat{\psi}_1^*-
   \nabla \hat{\psi}_2^*
     +\nabla \hat{\psi}_3^*) \cdot
   \ln \Psi
     \nabla \D \Psi}
     \\ & \tex{
      - c_2 \int\limits_{\ren} (\nabla \hat{\psi}_1^*+
   \nabla \hat{\psi}_2^*-\nabla \hat{\psi}_3^*) \cdot
   \ln \Psi \nabla \D \Psi
   } \\ & \tex{
   + \int\limits_{\ren} (
   \nabla \hat{\psi}_2^*-\nabla \hat{\psi}_3^*) \cdot
   \ln \Psi- \frac{\a_2}{4} \int\limits_{\ren} (\hat{\psi}_2^*-\hat{\psi}_3^*)
   y \cdot\nabla \hat{\psi}_1 
   }
\end{align*}
  and
\begin{align*}
    \tex{
    \o_2 (c_2,c_3)} & \tex{
    :=c_3 \int\limits_{\ren} \nabla \hat{\psi}_2^* \cdot
   \ln \Psi \nabla \D \Psi - c_2 \int\limits_{\ren} \nabla \hat{\psi}_3^* \cdot
   \ln \Psi \nabla \D \Psi
   }
\end{align*}
are the perturbations of the quadratic polynomials
 $$
 \mf{F}_i(c_2,c_3) :=
A_i c_2^2+ B_i c_3^2 + C_i c_2 +D_i c_3+ E_i c_2 c_3,
\quad \hbox{with} \quad i=1,2.
 $$
  The coefficients
of those quadratic expressions are given by
\begin{align*}
    &
     \tex{
     A_1:=  -\frac{\a_2}{4}\, \int\limits_{\ren} (\hat{\psi}_1^*+
    \hat{\psi}_2^*-\hat{\psi}_3^*) y
    \cdot(\nabla  \hat{\psi}_1-\nabla \hat{\psi}_2)  ,
 }
    \\ &
    \tex{
  B_1:= \frac{\a_2}{4}\, \int\limits_{\ren} (\hat{\psi}_1^*-
    \hat{\psi}_2^*+\hat{\psi}_3^*) y
    \cdot( \nabla \hat{\psi}_1 - \nabla\hat{\psi}_3) ,
    }
     \\ &
  C_1:=
 \tex{
   \frac{\a_2}{4}\,
    [\int\limits_{\ren} (\hat{\psi}_2^*-\hat{\psi}_3^*) y \cdot\nabla
    (2\hat{\psi}_1 -\hat{\psi}_2)- \int\limits_{\ren}
    \hat{\psi}_1 y \cdot\nabla  \hat{\psi}_1] ,
    }
     \\ &
  D_1:=
  \tex{
   \frac{\a_2}{4}\, [\int\limits_{\ren} (\hat{\psi}_2^*-\hat{\psi}_3^*) y \cdot\nabla
    (2\hat{\psi}_1 -\hat{\psi}_3)- \int\limits_{\ren}
    \hat{\psi}_1 y \cdot\nabla  \hat{\psi}_1] ,
   }
   \\ &
    \tex{
   E_1:=
    \frac{\a_2}{4}[\, \int\limits_{\ren} \hat{\psi}_1^* y
    \cdot( \nabla  \hat{\psi}_3-\nabla \hat{\psi}_2) - \, \int\limits_{\ren}
    (\nabla \hat{\psi}_2^* - \nabla\hat{\psi}_3^*) y \cdot (2\nabla \hat{\psi}_1 -
    \nabla \hat{\psi}_2-\nabla \hat{\psi}_3)],
   }
\end{align*}
\begin{align*}
    &
     \tex{
     A_2:=  -\frac{\a_2}{4}\, \int\limits_{\ren} \hat{\psi}_3^* y
    \cdot( \nabla  \hat{\psi}_1-\nabla \hat{\psi}_2),
 }
    \\ &
    \tex{
  B_2:= \frac{\a_2}{4}\, \int\limits_{\ren} \hat{\psi}_2^* y
    \cdot(
    (\nabla \hat{\psi}_1 - \nabla\hat{\psi}_3) ,
    }
     \\ &
  C_2:=
 \tex{
   \frac{\a_2}{4}\,
    \int\limits_{\ren} \hat{\psi}_3^* y \cdot\nabla  \hat{\psi}_1,
    }
     \\ &
  D_2:=
  \tex{
   -  \frac{\a_2}{4}\,
    \int\limits_{\ren} \hat{\psi}_2^* y \cdot\nabla  \hat{\psi}_1,
   }
   \\ &
    \tex{
   E_2:=
    \frac{\a_2}{4}\, \int\limits_{\ren} \hat{\psi}_2^* y
    \cdot (\nabla  \hat{\psi}_1-\nabla \hat{\psi}_2) -
    \hat{\psi}_3^* y
    \cdot (\nabla \hat{\psi}_1 - \nabla\hat{\psi}_3).
   }
\end{align*}

Therefore,  using the conic classification to solve \eqref{br64},
we will have the number of solutions through the intersection of
two conics. Then, depending on the type of conic, we shall always
obtain one to four possible solutions for our system. Hence,
somehow, the number of solutions depends on the coefficients we
have for the system and, at the same time, on the eigenfunctions
that generate the subspace $\ker\big({\bf B}+\frac{k}{4}\big)$.
Thus, we have the following conditions, which will provide us with
the conic section of each equation of the system \eqref{br64}:
\begin{enumerate}
\item[(i)] If $B_i^2 - 4A_iE_i < 0$, the equation represents an {\em ellipse},
unless the conic is degenerate, for example $c_2^2 + c_3^2 + k = 0$
for some positive constant k. So, if $A_i=B_j$ and $E_i=0$,
the equation represents a {\em circle};
\item[(ii)] If $B_i^2 - 4A_iE_i = 0$, the equation represents a {\em parabola};
\item[(iii)] If $B_i^2 - 4A_iE_i > 0$, the equation represents a {\em hyperbola}.
If we also have $A_i + E_i = 0$ the equation represents a
hyperbola (a rectangular hyperbola).
\end{enumerate}
Consequently, the zeros of the system \eqref{br64} and, hence, of
the system \eqref{br61}, adding the ``normalizing" constraint
\eqref{br15}, are ascertained by the intersection of those two
conics in \eqref{br64} providing us with the number of possible
$n$-branches between one and {four}. Note that in case  those
conics are two circles we only have two intersection points at
most. Moreover, due to the dimension of the eigenspaces it looks
like in this  case that we have four possible intersection points
two of them will coincide. However, the justification for this is
far from clear.

Moreover, as was done for the
previous case when $|\b|=1$, we need to control the oscillations
of the perturbation functions in order to maintain the number of
solutions. Therefore, imposing that
\begin{equation*}
    \left\| \o_i (c_2,c_3) \right\|_{L^\infty} \leq \mf{F}_i(c_2^*,c_3^*),
    \quad \hbox{with}\quad i=1,2,
\end{equation*}
we ascertain that the number of solutions must be between one and four.
This again gives us an idea of the difficulty of more general
multiplicity results.


\subsection{Further comments on mathematical justification of existence}
 \label{S6}


We return to the self-similar nonlinear eigenvalue  problem
\eqref{sf5}, associated with \eqref{i1},
which can be written in the  form
\begin{equation}
\label{cp7}
 \tex{
 {\mathcal L}(\a,n)f + {\mathcal N}(n,f)=0  \whereA
    {\mathcal N}(n,f):=\nabla \cdot ((1-|f|^n)\nabla \D f)\,.
    }
\end{equation}
As we have seen, the main difficulty in justifying the
$n$-branching behaviour concerns the distribution and
``transversal topology" of zero surfaces of solutions close to
finite interface hyper-surfaces.

\par
 Recall that,  as in classic nonlinear operator theory \cite{Deim, KZ,
 VainbergTr}, our  analysis above always assumed that we actually
 dealt with and performed computations for the
 integral equation:
\begin{equation}
\label{cp8}
\tex{
    f= -\mathcal{L}^{-1}(\a,n) \mathcal{N}(n,f) \equiv \mathcal{G}(n,f),\quad
    \mathcal{L}(\a,n):= -\D^2 +\frac{1-\a n}{4}\, y \cdot \nabla +\a I ,
    }
\end{equation}
where $\mathcal{L}(\a,n)$ is invertible in $L^2_\rho$ (this is
directly checked via Section \ref{S3}) and, hence compact, for a
fixed $\a$, and $f \in C_0(\ren)$ for small $n>0$. This confirms
that the zeros of the function $\mathcal{F}(n,f)$ are fixed points
of the map $\mathcal{G}(n,f)$. Note again that \eqref{cp8} is an
eigenvalue problem, where admissible real values of $\a$ are
supposed to be defined together with its solvability.
 This makes existence/multiplicity questions for \eqref{cp8}
 extremely difficult.
\par

There are two cases of this problem. The first and simpler one
occurs when the eigenvalue $\a$ is determined {\em a priori},
e.g., in the case $k=0$, where $\a_0(0)= \frac N{4}$ denoted as
$\a_0(0)=\a_0$, and where, for $n>0$, the first nonlinear
eigenvalue is given explicitly (see (\ref{alb1})):
 $$
 \tex{
 \a_0(n)= \frac N{4+Nn}.
 }
 $$
\par
Then \eqref{cp8} with $\a=\a_0(n)$ for $n>0$ becomes a standard
nonlinear integral equation with, however,  a quite curious and
hard-to-detect functional setting. Indeed, the right-hand side in
\eqref{cp8}, where the nonlinearity is not in a fully divergent
form, assumes the extra regularity at least such as
 \begin{equation}
 \label{ppp1}
 f \in H_\rho^3.
 \end{equation}
In view of the known good properties of the compact resolvent
$(\mathcal{L}-\l I)^{-1}$, it is clear that the action of the
inverse one $\mathcal{L}^{-1}$ is sufficient to restore the
regularity, since locally in $\ren$ this acts like $\D^{-2}$.
Therefore, it is plausible that
 \begin{equation}
 \label{ppp2}
\mathcal{G}: H^3_\rho \to H^3_\rho,
 \end{equation}
 and it is not difficult to get an {\em a priori} bound at least
 for small enough $f$'s. The accompanying  analysis as $y \to \infty$ (due to the unbounded domain)
  assumes no
 novelties or special difficulties and is standard for such
 weighted $L^2$ and Sobolev spaces.

  Therefore, application of Schauder's Fixed Point Theorem
 (see e.g., \cite[p.~90]{Berger})
   to
  \eqref{cp8} is a most powerful tool to imply existence of a
  solution, and moreover a continuous curve of fixed points $\Gamma_n=\{f,
  \,\, n>0\,\,\mbox{small}\}$. By scaling invariance of the
  similarity equation, we are obliged to impose the normalization
  condition, say,
   \begin{equation}
   \label{ppp3}
   f(0)= \delta_0>0 \quad \mbox{sufficiently small}.
    \end{equation}

 Uniqueness remains a completely open problem. However, studying the
 behaviour of the solution curve $\Gamma_n$ as $n \to 0$ and
 applying (under suitable hypothesis) the branching techniques
 developed above,
  we may conclude that any such continuous
 curve must be originated at a properly scaled eigenfunction
 $\psi_0=F$, so that such a curve is unique due to well-posedness
 of all the asymptotic expansions.

 A possibility of extension of $\Gamma_n$ for larger values of
 $n>0$ represents an essentially more difficult nonlocal open problem. Indeed, via
 compactness of linear operators involved in \eqref{cp8}, it is
 easy to expect that such a curve can end up at a bifurcation
 point only (unless blows up). However, nonexistence of turning
 saddle-node points at some $n_*>0$ (meaning that the $n$-branch is nonexistent for
 some $n>n_*$) is not that easy to rule
 out. Moreover, such turning points with thin film
 operators involved
 are actually possible, \cite{GalPetII}.

After establishing existence of such solutions for small $n>0$, we
face the next problem on their asymptotic properties including the
fact that these are compactly supported. On a qualitative level,
these questions were discussed in \cite{EGK1}.

In the case of higher-order nonlinear eigenfunctions of
\eqref{cp8} for $k \ge 1$ including the dipole case $k=1$, the
parameter $\a$ becomes an eigenvalue that is essentially involved
into the problem setting. This assumes to consider the equation
\eqref{cp8} in the extended space
 \begin{equation}
 \label{ppp4}
 (f,\a) \in X= H^3_\rho \times \{\a \in \re\} \,\,\,
 \mbox{and}
 \,\,\,
 \mathcal{G}: X \to X,
  \end{equation}
  where  proving the latter  mapping for some compact subsets becomes a
  hard
  open problem. Note that here even the necessary convexity issue
  for applying Schauder's Theorem
  can be  hard.
   We still do not know whether the representation such
  as \eqref{ppp4} may lead to any rigorous treatment of the
  nonlinear eigenvalue problem \eqref{cp8} for $k \ge 1$.

\section{General  Cauchy problem: a homotopic approach}
 \label{SHom7}

\subsection{Key concepts to justify: a first discussion}

 We now discuss some related properties of more general solutions of the CP for
the TFE \eqref{i1} using a homotopic approach when the parameter
$n$ approaches zero.
As shown in Section \ref{S3}, we already know the similarity
expression for the solutions of the ``limiting" bi-harmonic
equation \eqref{s1}. This fundamental solution has also a
self-similar structure thanks to the scaling invariance and the
uniqueness of the fundamental solution of the equation \eqref{s1},
denoted by \eqref{s3}.
\par
The idea is to perform a homotopic approach from \eqref{i1} to
\eqref{s1} in order to reveal  important (and still obscure in
general)  properties of the Cauchy problem. The reason is that the
bi-harmonic equation (\ref{s1}) i.e., \eqref{i1} when $n=0$, with
the same initial data, admits the unique classic solution given by
the convolution (\ref{s2}), where $b(x,t)$ is the fundamental
solution (\ref{s3}) of the operator $\frac{\p}{\p t} + \D^2$,
 and the oscillatory rescaled  kernel $F(y)$ is the unique solution of the
problem \eqref{s4}. Hence, we expect that the knowledge of the
solutions of \eqref{s1} can be extended to \eqref{i1} at least for
sufficiently small $n>0$.
 In other words, we claim that the ``fundamental" solutions for $n=0$
 and small $n>0$ exhibit several similar properties, excluding, on
 the other hand, some others such as the compact support one for $n>0$.
 In addition, the homotopic path $n \to 0^+$ can be used for a
 proper definition of the solutions of the Cauchy problem for the
 TFE--4 (\ref{i1}).

\par
 Thus, we  assume that $n>0$ is sufficiently small. We define
 some  ``homotopic classes" of degenerate parabolic PDE's
saying that the TFE \eqref{i1} is homotopic to the linear PDE
\eqref{s1} if there exists a family of uniformly parabolic
equations (a {\em homotopic deformation}) with coefficient
$\phi_\e(u)$ analytic in both variables $u\in\re$ and $\e\in (0,
1]$,
\begin{equation}
\label{cp4}
    u_{t} = -\nabla \cdot (\phi_\e(u)\nabla \D u)\,,
\end{equation}
such that $\phi_1(u)=1$ and
\begin{equation}
\label{cp5}
    \phi_\e(u)\rightarrow |u|^n\quad \hbox{as}\quad \e \rightarrow 0 \quad
    \hbox{uniformly on compact subsets}\,.
\end{equation}
We should point out that such a limit for {\em nonnegative} and
not changing sign solutions,  with various
 non-analytic (and non-smooth) regularizations has been widely used before in TFE--FBP theory as
 a key foundation; cf. \cite{BBP}, \cite{BF1}, and
\cite{BGK}.
\par
A possible homotopic path can be
\begin{equation*}
    \phi_\e(u):=\e^n+(1-\e)(\e^2+u^2)^{\frac n2}, \quad \e \in (0, 1].
\end{equation*}
For any $\e\in (0, 1]$, denote by $u_\e(x, t)$ the unique solution
of the CP for the regularized nondegenerate equation \eqref{cp4}
with same data $u_0$. By classic parabolic theory, $u_\e$ is
continuous (and analytic) in $\e\in (0, 1]$ in any natural
functional topology. The main problem is the behaviour as $\e
\rightarrow 0$, where the regularized PDE loses its uniform
parabolicity. For second-order parabolic equations obeying the
Maximum Principle, such regularization-continuity approaches are
typical for constructing unique solutions with singularities
(finite time blow-up, extinction, finite interfaces, etc.); see
\cite{GalGeom} as a source of key references and basic results.
However, for higher-order degenerate parabolic flows admitting
strongly oscillatory solutions of changing sign, such a
homotopy-continuity approach generates a number of difficult
problems. In fact, despite the fact that the passage to the limit
as $\e \to 0$ looks like a reasonable way to define a proper
solution of the TFE, we expect that there are always special
classes of compactly supported initial data, for which such a
limit is non-existent and, moreover, there are many partial
limits, thus defining a variety of different solutions (meaning
nonuniqueness), as we show below.

\subsection{Preliminary estimates}

 To ascertain such a limit for
\eqref{cp4} when $\e \rightarrow 0$, we firstly obtain some
estimations for its regularized solutions $\{u_\e(x, t)\}$. Here,
by $\O$ we denote either $\ren$, or, equivalently, the bounded
domain $\G_0\cap\{t\}$, i.e., the section of the support.

\begin{proposition}
\label{Pr cp1} Let $u_\e(x, t)$ be the unique global solution of
the CP for the regularized nondegenerate equation \eqref{cp4} with
the initial data $u_0$. Then, for any $t\in[0,T]$, the following
is satisfied:
\begin{enumerate}
\item[(i)] $u_\e(\cdot,t) \in H_0^1(\O)$;
\item[(ii)] $u_\e(\cdot,t) \in L^p(\O)$, with $p=1,2,\infty$; and
\item[(iii)] $h_\e \in L^2(\O \times [0,T])$, with
$h_\e := \phi_\e(u_\e) \nabla \D u_\e$.
\end{enumerate}
\end{proposition}

  \noi{\em Proof.}
Firstly, multiplying \eqref{cp4} by $\D u_\e$, integrating in $\O
\times [0,t]$ for any $t\in[0,T]$, and applying the formula of
integration by parts yield
\begin{equation}
\label{cp6}
 \tex{
    \frac{1}{2}\, \int\limits_\O |\nabla u_\e(x,t)|^2 +
    \int\limits_0^t \int\limits_\O \phi_\e(u)|\nabla \D u_\e|^2 =
    \frac{1}{2}\, \int\limits_\O |\nabla u_\e(x,0)|^2,
    }
\end{equation}
thanks to the boundary conditions \eqref{i3}. Note that
\begin{equation*}
 \begin{matrix}
      \int\limits_\O [|\nabla u_\e(x,t+h)|^2 - |\nabla u_\e(x,t)|^2]
 \ssk\ssk\\
      =\,
      -\int\limits_\O [\D u_\e(x,t+h)+ \D u_\e(x,t)][u_\e(x,t+h)-u_\e(x,t)].
      \end{matrix}
\end{equation*}
Then, dividing that equality by $h$, passing to the limit as
$h\downarrow 0$, and integrating between 0 and any $t\in [0,T]$,
we find that
\begin{equation*}
 \tex{
    \int\limits_0^t \int\limits_\O \D u_\e u_t =\frac{1}{2}\, \int\limits_\O |\nabla u_\e(x,t)|^2 -
    \frac{1}{2}\, \int\limits_\O |\nabla u_\e(x,0)|^2,
    }
\end{equation*}
which provides us with the necessary expression to obtain
\eqref{cp6}. Thus, from \eqref{cp6}, we have that (in fact, this
is true from the beginning for classic $C^\infty$-smooth solutions
of (\ref{cp4}), but we will need those manipulations in what
follows)
\begin{equation}
\label{cp33}
 \tex{
    \int\limits_\O |\nabla u_\e(x,t)|^2 \leq K
    \quad \hbox{and} \quad
    \int\limits_0^t \int\limits_\O \phi_\e(u)|\nabla \D u_\e|^2 \leq
    K,
    }
\end{equation}
since both terms of the left-hand side in \eqref{cp6} are always
positive and the right-hand side is bounded by \eqref{i4}, for
some positive constant $K>0$ that is independent of $\e$. Then,
\begin{equation*}
    \nabla u_\e(\cdot,t) \in L^2(\O) \quad \hbox{for any}\quad
    t\in[0,T].
\end{equation*}
Moreover, by Poincar\'e's inequality, $u_\e(\cdot,t) \in L^2(\O)$,
and   hence,
\begin{equation}
\label{cp10}
    u_\e(\cdot,t) \in H_0^1(\O) \quad \hbox{for any}\quad
    t\in[0,T].
\end{equation}
In fact, we may assume that
 \begin{equation}
 \label{Linf1}
  u_\e(\cdot,t) \in L^\infty(\O)
 \quad \mbox{for all} \quad t \in(0,T).
  \end{equation}
  Indeed, for $N=1$, this follows from (\ref{cp33}) by Sobolev's
  embedding. For $N \ge 2$, this is a natural assumption
inherited from the smooth analytic parabolic flow (\ref{cp4}),
though its full proof sometimes can be a difficult issue; we refer
to scaling and other techniques that may be convenient,
\cite{GMPKS}. From the conservation of mass assumption, we can
also assure that $u_\e(\cdot,t) \in L^1(\O)$. Also, expression
\eqref{cp6} combined with $u_\e(\cdot,t) \in L^1(\O)$ provides us
with the estimate
\begin{equation}
\label{cp12}
    h_\e \in L^2(\O \times [0,T])
    \quad(h_\e = \phi_\e(u_\e) \nabla \D u_\e).
\end{equation}
Indeed, from \eqref{cp33} we find that
\begin{equation*}
 \tex{
    \int\limits_0^t \int\limits_\O
    [\e^n+(1-\e)(\e^2+u^2)^{\frac n2}]|\nabla \D u_\e|^2 \leq
    K, \quad \mbox{so that}
    }
\end{equation*}
\begin{equation}
\label{cp52}
 \tex{
    \e^n\int\limits_0^t\int\limits_\O |\nabla \D u_\e|^2 \leq K
    \quad \hbox{and} \quad
    \int\limits_0^t \int\limits_\O
    (\e^2+u^2)^{\frac n2}|\nabla \D u_\e|^2 \leq
    K,
    }
\end{equation}
\noi since $\e \in (0,1)$ with a constant $K>0$ independent of
$\e$. Now, using H\"older's inequality,
\begin{equation*}
 \tex{
    \int\limits_0^t \int\limits_\O |h_\e|^2 \leq
    2\e^{2n} \int\limits_0^t \int\limits_\O
    |\nabla \D u_\e|^2 +2 \int\limits_0^t \int\limits_\O
    (\e^2+u^2)^{\frac n2}(\e^2+u^2)^{\frac n2}|\nabla \D u_\e|^2,
    }
\end{equation*}
by \eqref{cp52} and (\ref{Linf1}) (note also that  $u_\e(\cdot,t)
\in L^1(\O)$), we obtain \eqref{cp12}.
$\qed$

\ssk

Furthermore, the following estimates are also ascertained:

\begin{lemma}
\label{Le cp2} Let $u_\e(x, t)$ the unique global solution of the
CP for the regularized uniformly parabolic  equation \eqref{cp4}
with the initial data $u_0$. Then, there exists some positive
constant $K>0$ such that, for $x_1,x_2 \in \O$, independently of
$\e$ and t,
\begin{equation*}
     |u_\e(x_1,t)-u_\e(x_2,t)| \leq K |x_1-x_2|^{\frac 12}
      \quad \mbox{for odd $N$, \, and}
\end{equation*}
\begin{equation*}
     |u_\e(x_1,t)-u_\e(x_2,t)| \leq K |x_1-x_2|
     \quad \mbox{for even $N$}.
\end{equation*}

\end{lemma}

\noi{\em Proof.}
 Thanks to Sobolev's inequality, we have that, for the odd dimension $N\geq
 1$,
\begin{equation*}
    \big\| u_\e(\cdot,t)\big\|_{C^{0,\frac 12}(\O)}\leq C
    \big\| u_\e(\cdot,t)\big\|_{H_0^{N-j_N}(\O)}, \quad \hbox{with}
    \quad j_1=0,j_3=1,j_5=2,\cdots,
\end{equation*}
for a positive constant $C>0$. On the other hand, when the
dimension $N\geq 1$ is even,
\begin{equation*}
    \big\| u_\e(\cdot,t)\big\|_{C^{0,1}(\O)}\leq C
    \big\| u_\e(\cdot,t)\big\|_{H_0^{N-j_N}(\O)}, \quad \hbox{with}
    \quad j_2=0,j_4=1,j_6=2,\cdots,
\end{equation*}
for some positive constant $C>0$. Hence, since $C_0^\infty(\O)$ is
dense in $W_0^{k,2}(\O) = H_0^{k}(\O)$ for $k\geq 1$, by the
analytic smoothness  of the solutions of the uniformly parabolic
PDE \eqref{cp4} with analytic coefficients, for some positive
constant $K$ independently of $\e$ and $t$,
\begin{equation*}
     |u_\e(x_1,t)-u_\e(x_2,t)| \leq K |x_1-x_2|^{\frac 12}
   \quad \mbox{for odd $N$, \, and}
\end{equation*}
\begin{equation*}
     |u_\e(x_1,t)-u_\e(x_2,t)| \leq K |x_1-x_2|
     \quad \mbox{for even $N$.}  \qed
\end{equation*}

\ssk

 Moreover, as was noted in
 \cite{BBP, BF1} and \cite{BGK}, $u_\e(x,t)$ is also H\"{o}lder continuous in time
with exponent $\frac 1 8$, i.e., $u_\e(\cdot,t) \in C^{0,\frac
18}([0,T])$ for the one dimensional case. However, we provide a
new version for the $N$-dimensional case.

\begin{lemma}
\label{Le cp3} Let $u_\e(x, t)$ be the unique global  solution of
the CP for the regularized  equation \eqref{cp4} with the initial
data $u_0$. Then, there exists a positive constant $\tau>0$ such
that
\begin{equation*}
     |u_\e(x,t_2)-u_\e(x,t_1)| \leq \tau |t_2-t_1|^{\frac N{4(N+1)}},
\end{equation*}
if the dimension N is odd, and
\begin{equation*}
     |u_\e(x,t_2)-u_\e(x,t_1)| \leq \tau |t_2-t_1|^{\frac N{3N+2}},
\end{equation*}
if the dimension N is even,
for any $t_1,t_2 \in [0,T]$, independently of $\e$ and $x$.
\end{lemma}

\noi{\em Proof.} First, consider a non-negative cut-off function
$\var \in C_0^\infty (\re^N)$ such that ${\rm supp}\, (\var)
\subset \O$ and $\int_\O \var =1$. Subsequently, multiplying
\eqref{cp4} by a test function $\var_\g \in C_0^\infty(\re^N)$,
where
 $$
  \tex{
 \var_\g (x) := \frac 1 \g\,
\, \var\big(\frac {x-x_0}{\g^{1/N} }\big),
 }
 $$
with some $x_0 \in \O$ and  a constant $\g>0$ to be properly chosen later on,
 integrating over  $\O \times
[t_1,t_2]$ and applying the formula of integration by parts, we
find that
\begin{equation}
\label{cp50}
 \begin{matrix}
 \tex{
    -\int\limits_{t_1}^{t_2} \int\limits_\O \var_{\g} u_{\e,t} -
    \int\limits_{t_1}^{t_2} \int\limits_\O \nabla \var_\g \cdot (\phi_\e(u_\e)\nabla \D u_\e)
    =0,}
    \ssk \\
    \mbox{where}\quad -\int\limits_{t_1}^{t_2} \int\limits_\O \var_{\g} u_{\e,t}
    \equiv \int\limits_\O \var_\g (u_\e(t_2)-u_\e(t_1)).
    \end{matrix}
\end{equation}
On the other hand, by the present  choice of $\var_\g$, we know
that $\int_\O \var_\g =1$, since the Jacobian of $\frac
{x-x_0}{\g^{1/N} }$ is $\frac{1}{\g}$. Then, we have that
\begin{align*}
  &  u_\e(x_0,t_2)-u_\e(x_0,t_1) \equiv
  \tex{
  \int\limits_\O
    \var_\g(x)(u_\e(x_0,t_2)-u_\e(x_0,t_1))
    }
    \\
   & \equiv
   \tex{
   \int\limits_\O
   }
    \var_\g(x)(u_\e(x_0,t_2)-u_\e(x,t_2)
    +u_\e(x,t_2)-u_\e(x,t_1)+u_\e(x,t_1)-u_\e(x_0,t_1))
   \\
     & \le
 \tex{
    \int\limits_\O \var_\g(x)
    |u_\e(x,t_2)-u_\e(x_0,t_2)|
    }
     +
 \tex{
      \big|\int\limits_\O \var_\g(x)
    (u_\e(x,t_2)-u_\e(x,t_1))\big|
    }
    \\ &
      +
 \tex{
     \int\limits_\O \var_\g(x))
    |u_\e(x_0,t_1)-u_\e(x,t_1)|.
    }
\end{align*}
Owing to Lemma\;\ref{Le cp2} and \eqref{cp50} on the time interval
$(t_1,t_2)$, taking into account in the last inequality,
 we obtain
\begin{equation*}
    |u_\e(x_0,t_2)-u_\e(x_0,t_1)| \leq 2 K |x-x_0|^{\frac 12}  +
     \tex{
    \Big|\int\limits_{t_1}^{t_2} \int\limits_\O \nabla \var_\g \cdot (\phi_\e(u)\nabla \D u_\e)\Big|, \quad \hbox{with} \quad N\;\; \hbox{odd},
    }
\end{equation*}
\begin{equation*}
    |u_\e(x_0,t_2)-u_\e(x_0,t_1)| \leq 2 K |x-x_0|  +
     \tex{
    \Big|\int\limits_{t_1}^{t_2} \int\limits_\O \nabla \var_\g \cdot (\phi_\e(u)\nabla \D u_\e)\Big|,
     \quad \hbox{with} \quad N\;\; \hbox{even}.
    }
\end{equation*}
Moreover, by H\"{o}lder's inequality and the choice of $\var_\g$,
\begin{equation*}
 \tex{
   \Big| \int\limits_{t_1}^{t_2} \int\limits_\O \nabla \var_\g \cdot (\phi_\e(u)\nabla \D
   u_\e)\Big|
    \leq \Big(\int\limits_{t_1}^{t_2} \int\limits_\O
    |\phi_\e(u)\nabla \D u_\e|^2\Big)^{\frac 12}
    \Big(\int\limits_{t_1}^{t_2} \int\limits_\O |\nabla \var_\g|^2\Big)^{\frac 12},
    }
\end{equation*}
and, hence, thanks also to Proposition \ref{Pr cp1},
\begin{equation*}
 \tex{
   \Big| \int\limits_{t_1}^{t_2} \int\limits_\O \nabla \var_\g \cdot (\phi_\e(u)\nabla \D u_\e)
   \Big|
    \leq \tau \g^{-\frac{N+2}{2N}} |t_2-t_1|^{\frac 12}.
    }
\end{equation*}
Therefore, overall, interchanging $x$ and $x_0$ yields
\begin{equation*}
    |u_\e(x,t_2)-u_\e(x,t_1)| \leq 2 K |x-x_0|^{\frac 12}  +
    \tau \g^{-\frac{N+2}{2N}} |t_2-t_1|^{\frac 12},
    \quad \hbox{with} \quad N\;\; \hbox{odd},
\end{equation*}
\begin{equation*}
    |u_\e(x,t_2)-u_\e(x,t_1)| \leq 2 K |x-x_0|  +
    \tau \g^{-\frac{N+2}{2N}} |t_2-t_1|^{\frac 12},
    \quad \hbox{with} \quad N\;\; \hbox{even}.
\end{equation*}
Thus, taking $2 K < \tau$, $|x-x_0|< \g$, and $\g <
|t_2-t_1|^{\b}$, we obtain that
\begin{equation*}
    |u_\e(x,t_2)-u_\e(x,t_1)| \leq \tau |t_2-t_1|^{\frac \b2}  +
    \tau |t_2-t_1|^{-\frac{(N+2)\b}{2N}+\frac 12},
    \quad \hbox{with} \quad N\;\; \hbox{odd},
\end{equation*}
\begin{equation*}
    |u_\e(x,t_2)-u_\e(x,t_1)| \leq \tau |t_2-t_1|^\b  +
    \tau |t_2-t_1|^{-\frac{(N+2)\b}{2N}+\frac 12},
    \quad \hbox{with} \quad N\;\; \hbox{even}.
\end{equation*}
Consequently, taking $\b = \frac{N}{2(N+1)}$, if $N$ is odd, and
$\b = \frac{N}{3N+2}$ when $N$ is even, completes the proof.
 $\qed$

\subsection{Passing to the limit}

To conclude this section, we show existence of weak solutions for
the degenerate parabolic problem \eqref{i1} passing to the limit
as $\e$ goes to zero. However, we must admit from the beginning
that, from the analysis performed below, it is not possible to
assure which the limit will be (the solution of some CP or maybe
the solution of some FBP).
\par

By Proposition\;\ref{Pr cp1}, since, for bounded supports $\O$,
the embedding $H_0^1(\O) \hookrightarrow L^2(\O)$ is compact, we
can extract a convergent subsequence in $L^2(\O)$ as $\e\downarrow
0$ for the solutions of \eqref{cp4} labelled again $u_\e(x,t)$
such that
\begin{equation}
\label{cp11}
    \lim_{\e \rightarrow 0} \big\| u_\e(\cdot,t) - U(\cdot,t)\big\|_{L^2(\O)} =0.
\end{equation}
Consequently, the convergence of the non-degenerate solutions of
the problem \eqref{cp4} is strong in $L^2(\O)$.

Moreover, thanks to the H\"{o}lder continuity proved in
Lemmas\;\ref{Le cp2} and \ref{Le cp3}, we have a strong
convergence
   as $\e \downarrow 0$
in $C^{0,\frac{1}{2},\frac{N}{4(N+1)}}(\bar \O\times [0,T])$, when
$N$ is odd, and in $C^{0,1,\frac{N}{3N+2}}(\bar \O\times [0,T])$,
when $N$ is even. This is possible after applying the
Ascoli--Arzel\'a Theorem, since $\{u_\e\}$ is uniformly bounded
and equicontinuous in $\bar \O\times [0,T]$. Of course, these
estimates can imply other even stronger convergence results, which
are not treated below in detail.

Note that one difficulty we face is whether this limit depends on
the taken subsequence or not. In other words, this analysis does
not include
 any uniqueness result, which is expected to be a more difficult
 open problem for such nonlinear degenerate parabolic TFEs in
 non-fully divergence form and with non-monotone operators.
 However, the principal issue of the analytic regularization via
 \eqref{cp4} is that it is expected to lead to a smoother solution at the
 interface than those for the standard FBP. The difference is that
 the analytic regularized family $\{u_\e\}$, in addition to
 \eqref{i3}, is assumed to guarantee that, a.e. on the interface
 (assumed now sufficiently smooth),
  \be
  \label{int4}
  \tex{
   \frac{\p^2 u}{\p {\bf n}^2} =0.
   }
    \ee
In fact, proper oscillatory solutions of the CP are assumed to
exhibit even more regularity at smooth interfaces \cite{EGK2}:
 \be
 \label{int5}
\tex{
   \frac{\p^l u}{\p {\bf n}^l} =0, \quad \mbox{where} \quad
   l=\big[\frac 3n \big]-1.
   }
    \ee
Therefore, as $n \to 0^+$, the smoothness of such solutions at the
interfaces increases without bounds. Obviously, this is not the
case for the FBP (a ``positive obstacle" one) with standard
conditions \eqref{i3} and a usual quadratic (``parabolic") decay
at the interfaces.

\par

Thus, as above and customary, multiplying \eqref{cp4} by a test
function $\var \in C_0^\infty(\bar\O \times (0,T))$ and
integrating by parts in $\O \times [0,T]$
   gives
\begin{equation*}
 \tex{
    -\int\limits_0^T \int\limits_\O \var_t u_\e -
    \int\limits_0^T \int\limits_\O \nabla\var \cdot (\phi_\e(u)\nabla \D u_\e)
    =0.
    }
\end{equation*}
 Next, operating with this equality, we find that
\begin{equation}
\label{cp13}
 \tex{
    \int\limits_0^T \int\limits_\O \var_t u_\e + \e^n \int\limits_0^T \int\limits_\O
    \nabla\var \cdot \nabla \D u_\e
    +(1-\e)\int\limits_0^T \int\limits_\O \nabla\var \cdot (
    (\e^2+u^2)^{\frac n2}\nabla \D u_\e) =0.
    }
\end{equation}
Applying H\"{o}lder's inequality, it is clear from \eqref{cp6}
that there exists a subsequence labeled by $\{\e_k\}$ such that
the second term of \eqref{cp13} approximates zero as $\e_k
\downarrow 0$ for a sufficiently small, $n\approx0$, $n>0$,
 $$
 \tex{
    \Big|\e_k^n \int\limits_0^T \int\limits_\O \nabla\var \cdot \nabla \D u_{\e_k}\Big|
    }
 \tex{
    \leq \e_k\Big(\e_k^{2(n-1)} \int\limits_0^T \int\limits_\O (\nabla \D
     u_{\e_k})^2\Big)^{\frac 12}
    \Big( \int\limits_0^T \int\limits_\O |\nabla\var|^2\Big)^{\frac 12}
 }
    \\
     \leq K \e_k \downarrow 0,
 $$
as $\e_k \downarrow 0$, for some positive constant $K>0$.

Moreover, on the subset $\mathcal{P}:=\{(x,t)\in \O\times
[0,T]\,;\, |u(x,t)|> \d>0\}$, for any arbitrarily small $\d>0$, it
is clear that the limiting solution as $n \to 0$ is a weak
solution of the bi-harmonic equation \eqref{s1}. Indeed, by the
regularity of the uniformly parabolic equation \eqref{cp4} and the
uniformly H\"{o}lder continuity of its solutions proved in
Lemmas\;\ref{Le cp2} and \ref{Le cp3}, we obtain that $u_{\e,t}$,
$\nabla u_{\e}$, $\D u_{\e}$, $\nabla \D u_{\e,x}$, and $\D^2
u_\e$ converge uniformly on compact subsets of $\mc{P}$. In
general, it is not that difficult to see that, as $\e=\e_k \to 0$
(along the lines of classic results in \cite{BF1} and related
others), we obtain a weak solution of the TFE--4, i.e.,
\begin{equation}
\label{cp14}
 \tex{
    \int\limits_0^T \int\limits_\mathcal{P} \var_t U
    +\int\limits_0^T \int\limits_\mathcal{P} \nabla\var \cdot |U|^n \nabla \D U =0,
 }
\end{equation}
where $U(x,t)$ is the limit obtained through \eqref{cp11}. We
naturally assume that $\var \in C^\iy_0({\mathcal P})$.
\par

However, in the ``bad" subset $\{|u|\leq \d\}$, for any
sufficiently small $\d \geq 0$, we must take $\e>0$ sufficiently
small and depending on $\d$. Indeed, we take $\e$ such that
$0<\e\leq \d$. Thus, applying H\"{o}lder's inequality to the third
term in \eqref{cp13} over the subspace where $|u|\leq \d$, we have
that
\begin{align*}
 \tex{
    \Big|  \int\limits_0^T \int\limits_{\{|u|\leq \d\}}
 }
     &
      \tex{
       \nabla\var \cdot
    ( (1-\e) (\e^2+u_\e^2)^{\frac n2} \nabla \D u_\e)
    \Big|
    }
    \\ &
 \tex{
    \leq \Big(\int\limits_0^T \int\limits_{\{|u|\leq \d\}} |\nabla\var|^2 \Big)^{\frac 12}
     \Big(\int\limits_0^T \int\limits_{\{|u|\leq \d\}} (1-\e)^2 (\e^2+u_\e^2)^n
     |\nabla \D u_\e|^2\Big)^{\frac 12}.
     }
\end{align*}
Then, since $\var \in C_0^\infty(\bar\O \times (0,\infty))$ and
$\e\in (0,1)$, we get
\begin{equation*}
 \tex{
    \Big|  \int\limits_0^T \int\limits_{\{|u|\leq \d\}}
 }
      \tex{
      \nabla\var \cdot
    ( (1-\e) (\e^2+u_\e^2)^{\frac n2} \nabla \D u_\e) \Big|\leq
 }
 \tex{
     C \Big(\int\limits_0^T \int\limits_{\{|u|\leq \d\}} (1-\e) (\e^2+u_\e^2)^n
     |\nabla \D u_\e|^2\Big)^{\frac 12}
     }
\end{equation*}
for some positive constant $C>0$. Making use of the fact that
$|u|\leq \d$ and by \eqref{cp33}, we find that
\begin{align*}
 \tex{
    \Big|  \int\limits_0^T \int\limits_{\{|u|\leq \d\}}
 }
    &
 \tex{
    \nabla\var \cdot
    ( (1-\e) (\e^2+u_\e^2)^{\frac n2} \nabla \D u_\e) \Big|\leq
 }
    \\ &
     \tex{
     C \Big(\int\limits_0^T \int\limits_{\{|u|\leq \d\}} (1-\e) (\e^2+\d^2)^{\frac n2}
     (\e^2+u_\e^2)^{\frac n2}     |\nabla \D u_\e|^2\Big)^{\frac 12},
      }
\end{align*}
and, hence, using \eqref{cp52},
\begin{equation}
 \label{nn11}
 \tex{
    \Big|  \int\limits_0^T \int\limits_{\{|u|\leq \d\}} \nabla\var \cdot
    ( (1-\e) (\e^2+u_\e^2)^{\frac n2} \nabla \D u_\e) \Big|\leq
     C_1 \d^{\frac n2} \sim C_1 \e ^{\frac n2},
     }
\end{equation}
for some constant $C_1>0$
and taking $\e \sim \d$
sufficiently small.

\par

 Finally, the estimate \eqref{nn11} shows the actual rate of the
 limit as $n \to 0$,  together with $\e \to 0$, in the analytic
 approximating flow \eqref{cp4} to get in this limit weak (and
 hence classic by standard parabolic theory)
 solutions of the bi-harmonic equation \eqref{s1}. Namely, one has
 to have that
  \be
  \label{nn12}
   \tex{
  n=n(\e) \to 0 \,\,\, \mbox{such that} \quad \e^{\frac{n(\e)}2}
  \to 0 \,\,
  \Longrightarrow\,\, n(\e) \gg \frac 1{|\ln \e|} \,\, \mbox{as} \,\, \e \to
  0^+.
   }
   \ee

However, this is not the end of the problem: indeed, under the
condition \eqref{nn12} on the parameters, we definitely arrive at
the limit $\e,\, n(\e) \to 0$ to the weak solution of the
bi-harmonic equation written in the following ``mild" form:
\begin{equation}
\label{cp14N}
 \tex{
    \int
    _0^T \int
    _\O \var_t U
    +\int
    _0^T \int
    _\O \nabla\var \cdot \nabla \D U
    =0.
 }
\end{equation}
This is not a full definition of weak solutions, since it assumes
just a single integration by parts, so allows us also positive
solutions of the ``obstacle" FBP for \eqref{s1} with the
corresponding conditions \eqref{i3} (with $n=0$), which can be
constructed by ``singular" regularization as in \cite{BF1}.

Thus, unfortunately, our analysis still does not recognize the
desired difference between oscillatory solutions of the CP and
others (possibly positive ones) of the standard FBP and others
that can be posed for the TFE--4 \eqref{i1}. Nevertheless, this
first step in a homotopy analysis declares useful estimates and
bounds on the parameters of regularization such as \eqref{nn12},
which are absolutely necessary for passing to the limit to get
sign changing solutions of the linear bi-harmonic flow.

Then, the homotopy concept as a connection to the linear PDE
\eqref{cp4} can  describe the origin (at $n=0$) of the oscillatory
solutions of TFEs and hence establish a transition to the {\em
maximal regularity} of the solutions of  \eqref{i1}. Indeed,
inevitably, bearing in mind the oscillatory character of the
kernel $F(|y|)$ of the fundamental solution, the proper solutions
of the CP are going to be oscillatory near finite interfaces at
least for small $n > 0$.


\subsection{Final remark: towards a full homotopy approach}

 The main obstacle faced in obtaining a full exhibition of
such a homotopic approach, when $n \to 0^+$, is the study of
difficult pointwise limits near interfaces, where key
singularities (a kind of Riemann's problems) occur. In particular,
as a clue, let us mention a proper
 topology, in which we can define the homotopy. This
is clear when we transform \eqref{i1} (here, for simplicity,
avoiding $\e$-regularization as in \eqref{cp4}) into a
perturbation of the bi-harmonic equation \eqref{s1}:
\begin{equation}
\label{cp9}
    u_{t} = -\D^2 u +g_n(u), \quad \hbox{where} \quad
    g_n(u):=\nabla \cdot ((1-|u|^n)\nabla \D u)\,,
\end{equation}
which we write down as an integral equation using the compact
resolvent $(\BB-\l I)^{-1}$ and the semigroup from Section
\ref{S3}.
 We then deal with the integral equation
  \be
  \label{inteq1}
   \tex{
 u(t)= b(t)*u_0 + \int
  _0^t b(t-s)* g_n(u(s))\, {\mathrm d} s,
  }
  \ee
  where $b(t)$ is the fundamental solution \eqref{s3}.
 As usual, the integral form \eqref{inteq1} allows us to weaken the necessary treatment
of the perturbation $g_n(u)$, which is assumed to be small as $n
\to 0$. However, this does not rule out the principal difficulty
concerning such an unusual and very sensitive perturbation
$g_n(u)$.

Thus, the main open problem of the homotopy issues is as follows:
under which ``topology-functional-geometric" setting for admitted
solutions $u(x,t)$,
 \be
 \label{int44}
 \fbox{$
  \tex{
 \int
 _0^t b(t-s)* g_n(u(s))\, {\mathrm d} s \to 0 \quad \mbox{a.e.
 as}
 \quad n \to 0^+.
 }
 $}
 \ee
 In particular, it is not difficult to
see that, for sufficiently smooth functions $u$ with a finite
number of transversal zero surfaces uniformly in small $n \ge 0$,
we have that
\begin{equation*}
    g_n(u)\to 0 \quad \hbox{as} \quad n \to 0^+\,.
\end{equation*}
at least a.e., and in other natural (weighted) topologies
associated with the operator $\BB$ and/or others. This can  even
be true uniformly on compact subsets, if the differential
operators in $u$ are bounded on such special functions, whose
regularity and the ``geometric transversal structure" near zero
surfaces well-correspond to the desired maximal
regularity/structure that are generic for  the TFE solutions. But
such a detailed {\em a priori} information on solutions seems to
be excessive and not required. Taking  the full homotopic approach
as $\e\, \,, n \to 0$, we desperately need to understand the
structure of the zero surfaces, because of the oscillatory
behaviour of these solutions of changing sign.  Away from small
neighbourhoods of such zeros (zero curves or surfaces), there is
no any essential problem. The structural properties of zeros of
solutions of   the TFE--4 \eqref{i1} and the TFE--6 have been
discussed in \cite{EGK2, EGK4}, and the results therein inspire us
with a certain optimism concerning the correctness of the general
homotopy approach to the CP for the TFEs, though  difficulties are
far away from being properly settled in a general setting.


\end{document}